\documentclass[11pt]{amsart}
\usepackage{amssymb}

\usepackage{a4}
\usepackage[all]{xy}

\usepackage{color}

\newtheorem{defi}{Definition}[section]
\newtheorem{sinobservacion}[defi]{}
\newenvironment{sinob}{\begin{sinobservacion} \rm}{\end{sinobservacion} }
\newtheorem{coro}[defi]{Corollary}
\newtheorem{lem}[defi]{Lemma}
\newtheorem{nota}[defi]{Notation}

\newtheorem{rem}[defi]{Remark}
\newtheorem{prop}[defi]{Proposition}
\newtheorem{thm}[defi]{Theorem}
\newtheorem{ej}[defi]{Example}

\newcommand{\benu}{\begin{enumerate}}
\newcommand{\enu}{\end{enumerate}}

\date{\today}

\begin{document}
\title[On the degree in categories of complexes of fixed size]{On the degree in categories of complexes of fixed size}
\author[Chaio]{Claudia Chaio}
\address{Centro Marplatense de Investigaciones Matem\'aticas, Facultad de Ciencias Exactas y
Naturales, Funes 3350, Universidad Nacional de Mar del Plata, 7600
Mar del Plata, Argentina}
 {\email{claudia.chaio@gmail.com}
\author[Pratti]{Isabel Pratti}
\address{Centro Marplatense de Investigaciones Matem\'aticas, Facultad de Ciencias Exactas y
Naturales, Funes 3350, Universidad Nacional de Mar del Plata, 7600
Mar del Plata, Argentina} \email{nilprat@mdp.edu.ar}
\author[Souto Salorio]{Mar\'ia Jos\'e Souto Salorio}
\address{Facultade de Informatica, Campus de Elvi\~{n}a,
Universidade da Coru\~{n}a, CP 15071, A Coru\~{n}a, Espa\~{n}a}
{\email{maria.souto.salorio@udc.es}

\date{\today}

\thanks{The first and second authors thankfully acknowledge partial support from CONICET
and EXA558/14 from Universidad Nacional de Mar del Plata, Argentina. The third author thankfully acknowledge support
from Ministerio Espa\~{n}ol de Econom\'ia y Competitividad and FEDER (FF12014-51978-C2-2-R). The first
author is a researcher from CONICET}



\maketitle

\begin{abstract}
We consider $\Lambda$
an artin algebra and $n \geq 2$.  We study how to compute the left and right degrees of irreducible morphisms between complexes in a  generalized standard  Auslander-Reiten component of  ${\mathbf{C_n}({\rm proj}\, \Lambda)}$ with length.
We give conditions under which the kernel and the cokernel of irreducible morphisms between complexes in $\mathbf{C_n}({\rm proj}\, \Lambda)$ belong to such a category.
For a finite dimensional hereditary algebra  $H$ over an algebraically closed field, we determine when an irreducible morphism has finite left (or right) degree and we give a characterization, depending on the degrees of certain  irreducible morphisms, under which $\mathbf{C_n}({\rm  proj} \,H)$ is of finite type.
\end{abstract}

\section*{Introduction}

In the last years, the Auslander-Reiten theory has played an important role in the study of representations of artin algebras and additive categories.
The representation theory of finite dimensional algebras has been developed in the beginning of the seventies, with P. Gabriel results about the classifications of hereditary algebras 
which are representation-finite algebras over an algebraically closed field and, also with the  Auslander-Reiten theory, between others. The original ideas studied for categories of modules has been generalized to additive  categories with exact structure such as the categories of complexes.

In \cite{BSZ:04}, the authors studied the category of complexes of fixed size with projective $\Lambda$-modules entries which they denoted by $\mathbf{C_n}({\rm proj}\, \Lambda),$ with $\Lambda$  an artin algebra and $n \geq 2$. They proved that the mentioned category is a Krull-Schmidt category (an additive category in which each non-zero object is a finite direct sum of objects with local endomorphism algebra), has almost split sequences, is not abelian  and, moreover they observed that there is a strong relationship between the almost split sequences of complexes in $\mathbf{C_n}({\rm proj}\, \Lambda)$ and the Auslander-Reiten triangles in the bounded derived category  $D^b({\rm mod}\, \Lambda)$.

Recently, in \cite{CPS1}, the authors continued with the study of the category $\mathbf{C_n}({\rm proj}\, \Lambda)$ and they showed that the Auslander-Reiten quiver of $\mathbf{C_n}({\rm proj}\, \Lambda)$, denoted by $\Gamma_{\mathbf{C_n}({\rm  proj}\, \Lambda)}$, can be constructed by
the well-known "knitting algorithm" used to build the Auslander-Reiten quiver of ${\rm mod}\, \Lambda$, and furthermore they studied the behavior of sectional paths showing that their composition can be the zero morphism.

Since the category $\mathbf{C_n}({\rm proj}\, \Lambda)$ is not abelian, in this work we are interested to study when the kernel and the cokernel of irreducible morphisms between complexes in such a category belong to $\mathbf{C_n}({\rm proj}\, \Lambda)$, see Section 2. We recall that $f \colon X \rightarrow Y$, with $X,Y \in \mathbf{C_n}({\rm proj}\, \Lambda)$,
is an  irreducible morphism provided it does not split
and whenever $ f = gh$, then either $h$ is a section or $g$ is a
retraction, see $(1.2)$.

The notion of degree of an irreducible morphism was introduced by S. Liu for a Krull-Schmidt category in \cite{L:10}. This concept has shown to be a fundamental tool to solve many problems in representation theory of algebras.
For a category of finitely generated modules, we know that the kernel (cokernel) is fundamental when we want to determine the left (right) degree  of an irreducible epimorphism (monomorphism, respectively). In section 3, we study a similar problem for irreducible morphisms in $\mathbf{C_n}({\rm proj}\, \Lambda)$. More precisely, we study which complex is the domain of the morphism that determines the left degree for irreducible morphisms between complexes in generalized standard components of
$\Gamma_{\mathbf{C_n}({\rm proj}\, \Lambda)}$  with length.
We prove that, in case that the kernel of an irreducible morphism $f$  belongs to $\mathbf{C_n}({\rm proj}\, \Lambda)$ then this complex plays a  similar role as in ${\rm mod}\, \Lambda$, that is, if the inclusion  morphism $j:{\rm Ker}\, f \rightarrow X$ is such that $j \in \Re_{\mathbf{C_n}(\rm{proj}\,\Lambda)}^m({\rm Ker}\, f ,X)\backslash
\Re_{\mathbf{C_n}(\rm{proj}\,\Lambda)}^{m+1}({\rm Ker}\, f ,X)$ then the left degree of $f$ is $m$. Otherwise a right $\mathbf{C_n}(\rm{proj}\,\Lambda)$-approximation of the kernel has to be considered to determine the left degree. Precisely, in this second case we prove:
\vspace{.05in}

\noindent{\bf Theorem A.} {\it Let $\Lambda$ be an artin algebra and $\Gamma$  a generalized standard component
of $\Gamma_{\mathbf{C_n}(\rm{proj}\,\Lambda)}$ with length. Let $f:X \rightarrow Y$ be an irreducible morphism with $X \in \Gamma$,
$\varphi: W \rightarrow {\rm Ker}\,f$ a right
$\mathbf{C_n}(\rm{proj}\,\Lambda)$-approximation  of ${\rm Ker}\,f$ and $j:{\rm Ker}\, f \rightarrow X$ the inclusion morphism in $\mathbf{C_n}(\rm{mod}\,\Lambda)$. The following conditions are equivalent.
\begin{itemize}
\item [(a)] $d_l(f) = m$.
\item [(b)] $j \varphi \in \Re_{\mathbf{C_n}}^m(W,X)\backslash
\Re_{\mathbf{C_n}}^{m+1}(W,X)$. Moreover, there exists an indecomposable direct summand $Z$  of $W$ and a morphism $h \in   \Re_{\mathbf{C_n}(\rm{proj}\,\Lambda)}^m(Z,X)\backslash
\Re_{\mathbf{C_n}(\rm{proj}\,\Lambda)}^{m+1}(Z,X)$ such that $fh=0$.
 \end{itemize}}

Anyway we refer to $H$ in the text, we mean that $H$ is a hereditary algebra.

In \cite{CLT} the authors gave a  characterization in order to compute the degree of an irreducible morphism. To achieve such a result they used covering techniques. In section 4, we adapt such covering techniques in order to study the degree of irreducible morphisms between indecomposable complexes in a category $\mathbf{C_n}({\rm proj}\; H)$, for $H$ a finite dimensional hereditary algebra over an algebraically closed field. We obtain similar results as the ones for a category of finitely generated modules, but the proofs are not necessarily the same since the category we deal with has some particular properties, such as,  that the irreducible morphisms with injective indecomposable domain are not epimorphisms. We prove Theorem B.
\vspace{.05in}

\noindent{\bf Theorem B.} {\it Let $H$ be a hereditary finite dimensional $k$-algebra over an algebraically closed field and $f: X \rightarrow Y$ be an irreducible morphism in $\mathbf{C_n}({\rm proj} \,H)$ with
$X$ indecomposable. Let  $\Gamma$ be an Auslander-Reiten component  of
  $\Gamma_{\mathbf{C_n}({\rm proj} \,H)}$ containing
$X$ and $m \in \mathbb{N}$.
\begin{enumerate}
\item [(a)] If $d_l(f) =m$, then there exists $Z \in \Gamma$ and $h
\in \Re_{\mathbf{C_n}({\rm proj} \,H)}^{m} (Z,X) \backslash \Re_{\mathbf{C_n}({\rm proj} \,H)}^{m +1}(Z,X)$ such that $fh =0$.

\item [(b)] If $d_r(f) =m$, then there exists $Z \in \Gamma$ and $h
  \in \Re_{\mathbf{C_n}({\rm proj} \,H)}^{m} (Y,Z) \backslash \Re_{\mathbf{C_n}({\rm proj} \,H)}^{m
  +1}(Y,Z)$ such that $hf =0$.
\end{enumerate}}
\vspace{.05in}

Furthermore, we also
study equivalent conditions for a category $\mathbf{C_n}({\rm proj}\, H)$ to be representation-finite whenever $H$ is a hereditary algebra in terms of degrees of some particular irreducible morphisms.  Following the definitions and notations of  $(\ref{1.1})$ and $(\ref{T(S)})$ observe that if $P$ is a simple indecomposable  projective $H$-modulo then  $R_{n-1}(P)=T(P)$. Moreover, in case that $P$ is a projective $H$-module such that $\nu(P)$ is simple injective then  $L_{1}(P)=S(P)$. More precisely, we prove Theorem C.
\vspace{.05in}

\noindent{\bf Theorem C.} {\it Let $H$ be a hereditary finite dimensional $k$-algebra over an algebraically closed field. The following conditions are equivalents.
\begin{enumerate}
\item [(a)] The category $\mathbf{C_n}({\rm proj} \,H)$ is of finite type.
\item [ (b)]  There exists a simple projective $H$-module $P$ such that the irreducible morphism ${\rho}_{n-1}^n:T(P) {\rightarrow} J_{n-1}(P)$
         in $\mathbf{C_n}({\rm proj} \,H)$  has finite right degree.
\item [(c)] For all simple projective $H$-module $Q$,  the irreducible morphism
         ${\rho}_{n-1}^n:  T(Q) {\rightarrow} J_{n-1}(Q)$ in $\mathbf{C_n}({\rm proj} \,H)$ has finite right degree.
\item [(d)] There exists a  projective $H$-module $P$ with $\nu(P)$ a  simple injective $H$-module such that the irreducible morphism ${\lambda}^n_1: J_1(P){\rightarrow} S(P)$  in $\mathbf{C_n}({\rm proj} \,H)$ has finite left degree.
\item [(e)]  For all projective $H$-module $Q$ with $\nu(Q)$ a simple injective $H$-module, the irreducible morphism
        ${\lambda}^n_1: J_1(Q) {\rightarrow} S(Q)$  in
        $\mathbf{C_n}({\rm proj} \,H)$ has finite left degree.
\item [(f)] All irreducible morphisms $f:X \rightarrow Y$ with $X$ or $Y$ indecomposable complexes in $\mathbf{C_n}{(\rm proj} \,H)$ of type $(sec)$
        are such that $d_r(f) < \infty$.
\item [(g)] All irreducible morphisms $f:X \rightarrow Y$ with $X$ or $Y$ indecomposable complexes  in  $\mathbf{C_n}({\rm proj} \,H)$ of type $(ret)$ are such that $d_l(f)< \infty$.
\end{enumerate}}
\vspace{.05in}

In section 5, we end up with the study of the nilpotency index of the radical of a representation-finite category $\mathbf{C_n}({\rm proj} \,H)$ for $H$ a Dynkin algebra. We give a result depending on the nilpotency index of the radical of ${\rm mod}\,H$. Precisely, we prove Theorem D.
\vspace{.05in}

Let $H$ be a hereditary algebra of finite representation type, and $r_{_H}$ be the nilpotency index of the radical of ${\rm mod}\,H$. Let $\ell= {\rm max}\{\ell_{j,i}\}$ where $\ell_{j,i}$ is the number of arrows of a path from a source $i$ to a sink $j$ of $Q_{H}$.
\vspace{.05in}

\noindent{\bf Theorem D.} {\it Let $H$ be a hereditary algebra and $\mathbf{C_n}({\rm proj} \,H)$ be a category of finite type. The nilpotency index $r$ of the  radical of  $\mathbf{C_n}({\rm proj} \,H)$ is the following:
\begin{enumerate}
\item [a)] If  $n=2$  then $r=r_{_H} + \ell + 1$.
\item [b)] If  $n \geq 3$ then  $r= 2r_{_H} + \ell + 1$.
\end{enumerate}}

\section{Preliminaries}

Let $\Lambda$ be an artin algebra and ${\rm mod}\,\Lambda$ be the
category of finitely generated right $\Lambda$-modules. We denote
by ${\rm proj}\,\Lambda$  the full subcategory of ${\rm
mod}\,\Lambda$ consisting of all projective $\Lambda$-modules.

\begin{sinob} \label{1.1} We recall the definition of an object in the category of complexes $\mathbf{C}({\rm mod}\,\Lambda)$. A complex $X \in \mathbf{C}({\rm mod}\,\Lambda)$ is a collection, $(X^i,d_X^i)_{i \in \mathbb{Z}}$, with $X^i \in {\rm mod}\,\Lambda$ and $d^i_X \rightarrow X^{i+1}$ morphisms in ${\rm mod}\,\Lambda$ such that $d^{i+1}_X d_X^i=0$. If $X=(X^i,d_X^i)_{i \in \mathbb{Z}}$ and $Y=(Y^i,d_Y^i)_{i \in \mathbb{Z}}$ are two complexes, a morphism $f:X \rightarrow Y$  is a family of morphisms $\{f^i:X^i \rightarrow Y^i\}_{i \in \mathbb{Z}}$  in ${\rm mod}\,\Lambda$, such that  $f^{i+1}d_X^i=d_Y^if^i$, for all $i \in \mathbb{Z}$.

In the  category  $\mathbf{C}({\rm mod}\,\Lambda)$ we
consider the class $\mathcal{E}$  of composable morphisms
$X\stackrel{f}{\rightarrow }Y\stackrel{g}{\rightarrow }Z$
such that for all $i \in \mathbb{Z}$, we have an exact sequence $ 0\rightarrow
X^{i}\stackrel{f^{i}}{\rightarrow
}Y^{i}\stackrel{g^{i}}{\rightarrow }Z^{i}\rightarrow  0$ that splits.

For $n \in \mathbb{N}$, the category  ${\mathbf {C_n}}({\rm mod}\,\Lambda)$ is defined in \cite{BSZ:04}, as the full subcategory
of ${\mathbf C}({\rm mod}\,\Lambda)$ whose objects are those complexes
$X$ such that $X^{i}=0$ if $i\notin \{1, \dots , n\}$. Particularly, they studied the subcategory of ${\mathbf {C_n}}({\rm proj}\,\Lambda)$  whose objects  are those complexes $X$ in ${\mathbf {C_n}}({\rm mod}\,\Lambda)$ such that $X^i \in  {\rm proj}\,\Lambda$.
We denote by  $\mathcal{E}_{n},$  the class of
composable morphisms in $\mathbf{C_n}({\rm proj}\,\Lambda)$ which are in
$\mathcal{E}$, that is, for all $i=1,  \dots, n$, we have an exact sequence $ 0\rightarrow
X^{i}\stackrel{f^{i}}{\rightarrow
}Y^{i}\stackrel{g^{i}}{\rightarrow }Z^{i}\rightarrow  0$ that splits.
It is known that the category $(\mathbf{C_n}({\rm proj}\,\Lambda), \mathcal{E}_{n})$ is exact.
We denote by the expression
$X^1 {\rightarrow} X^{2} {\rightarrow} \cdots
{\rightarrow} X^n $ a complex in $\mathbf{C_n}({\rm proj} \,\Lambda)$.

Given an indecomposable module $P \in {\rm proj}\,{\Lambda}$, in \cite{BSZ:04} the following complexes are defined:
$J_{i}(P)=(J^{s},d^{s})_{s \in \mathbb{Z}}$ for $i \in \{1,
\cdots,n-1\}$, with $J^{s}=0$ if $s\neq i, s\neq i+1, J^{i}=J^{i+1}=P$,
$d^{i}=id_{P}$,
$T(P)=(Y^{i},d^{i})_{i\in \mathbb{Z}}$ with
 $Y^{i}=0$ for $i\neq n$, $Y^{n}=P$, $d^{i}=0$  for all  ${i\in \mathbb{Z}}$,
and $S(P)=(X^{i},d^{i})_{i\in \mathbb{Z}}$ with $X^{i}=0$ for $i\neq 1$, $X^{1}=P$,
$d^{i}=0$ for all ${i\in \mathbb{Z}}$.

By \cite[Corollary 3.3]{BSZ:04} the complexes $J_{i}(P)$ and $T(P)$ are the indecomposable projective objects in $\mathbf{C_n}({\rm proj}\, \Lambda)$ and, the complexes of the form $J_{i}(P)$ for $i=1, \dots, n-1$
and $S(P)$  are the indecomposable injective objects in $\mathbf{C_n}({\rm proj}\, \Lambda)$ where $P$ is an indecomposable module in ${\rm proj}\,{\Lambda}$.

In \cite[Proposition 3.6, 3.7]{BSZ:04}, the authors  proved that the category $\mathbf{C_n}({\rm proj} \, \Lambda)$ is exact with enough
projective and  injective objects.

Let $\mathbf{\overline{C_n}}({\rm  proj} \, \Lambda)$
(respectively, $\mathbf{\underline{C_n}}({\rm proj} \, \Lambda)$) be  the category whose objects are the ones of
$\mathbf{C_n}({\rm proj} \, \Lambda)$  and the morphisms are those in
$\mathbf{C_n}({\rm proj} \, \Lambda)$ that factor through injective complexes (respectively, projective complexes).

The category $\mathbf{C_n}({\rm proj} \, \Lambda)$  was first studied for  $n=2$ in \cite[Proposition 3.3]{B:04}, where an equivalence is established  between  ${\rm mod}\, \Lambda$ and  $\overline{\mathbf{C_2}}({\rm proj}\,\Lambda)$.

We recall that the shift $X[j]$ in $\mathbf{C_n}({\rm proj} \, \Lambda)$  is defined as follows; $(X[j])^t = X^{t+j}$ and $(d_{X[j]})^t = (-1)^j d^{t+j}$, for $j,t \in \mathbb{Z}$ whenever these expressions makes sense.

According to \cite[Remark 3.5]{BSZ:04} we consider $\mathbf{C}^{\leq n}({\rm proj}\,\Lambda)$ the full subcategory of $\mathbf{C}({\rm proj}\,\Lambda)$ whose objects are those $X \in \mathbf{C}({\rm proj}\,\Lambda)$ such that $X^i=0$ for $i>n$. Similarly, $\mathbf{C}^{\geq n}({\rm proj}\,\Lambda)$ is the full subcategory of $\mathbf{C}({\rm proj}\,\Lambda)$ whose objects are those $X \in \mathbf{C}({\rm proj}\,\Lambda)$ such that $X^i=0$ for $i < n$.

The functors $F: \mathbf{C}^{\leq n}({\rm proj}\,\Lambda)\rightarrow \mathbf{C_n}({\rm proj} \,\Lambda)$ and $G: \mathbf{C}^{\geq 1}({\rm proj}\,\Lambda) \rightarrow \mathbf{C_n}({\rm proj} \,\Lambda)$ are defined as follows: for $X \in \mathbf{C}^{\leq n}({\rm proj}\,\Lambda)$, we define  $F(X)=(W^i,d^i)_{i \in \mathbb{Z}}$ with $W^i=X^i$ for $i \geq 1$, $W^i=0$ for $i < 1$ and $d^i= d_X^i$, for $i \geq 1$. For $f:X \rightarrow Y$ a morphism in $\mathbf{C}^{\leq n}({\rm proj}\,\Lambda)$ we define $F(f)=(g^i)_{i \in \mathbb{Z}}$ with $g^i=f^i$ for $i = 1, \cdots, n$.
For $X \in \mathbf{C}^{\geq 1}({\rm proj}\,\Lambda)$ we define $G(X)=(W^i,d^i)_{i \in \mathbb{Z}}$ where $W^i= X^i$ for $i \leq n$, $W^i=0$ for $i > n$ and $d^i= d^i_X$ for $i \leq n-1$. If $f:X\rightarrow Y$ is a morphism in $\mathbf{C}^{\geq 1}({\rm proj}\,\Lambda)$, we define $G(f)=(g^i)_{i \in \mathbb{Z}}$ with $g^i= f^i$ for $i \leq n$.

We say that $X \in \mathbf{C}({\rm mod}\,\Lambda)$ has bounded cohomology if $H^i(X)=0$ for all but finitely many $i \in \mathbb{Z}$, where $H^i(X)= {\rm Ker}\,d_X^i/{\rm Im}\,d_X^{i-1}$. Denote by $\mathbf{C^{-,b}}({\rm mod}\,\Lambda)$ ($\mathbf{C^{+,b}}({\rm mod}\,\Lambda)$, respectively) the full subcategory of complexes bounded above with bounded cohomology (bounded below with bounded cohomology,respectively).

Recall that, $f$, $g \in {\rm Hom}_{\mathbf{C}({\rm mod}\,\Lambda)}(X,Y)$ are homotopic if there are morphisms $h^i:X^i \rightarrow Y^{i-1}$ such that $f^i - g^i= h^{i+1}d_X^i + d_Y^{i-1}h^i$, for all  $i \in \mathbb{Z}$.

Denote by $\mathbf{K}({\rm mod}\,\Lambda)$, $\mathbf{K^{-,b}}({\rm mod}\,\Lambda)$, $\mathbf{K^{+,b}}({\rm mod}\,\Lambda)$ and $\mathbf{K^b}({\rm mod}\,\Lambda)$ the homotopy categories of complexes introduced above. Similarly the above categories are defined for
 ${\rm proj}\,\Lambda$.

Following \cite[Notation 5.2.]{BSZ:04} we denote by $\mathfrak{L}_n$ the full subcategory of $\mathbf{K}({\rm proj}\,\Lambda)$ whose objects are those $X \in \mathbf{C}^{\leq n}({\rm proj}\,\Lambda)$ with $H^m(X)=0$ for $m \leq 1$.

In case $H$ is a hereditary algebra, for each $i=0, \cdots, n-2$ we have an embedding $\varphi_i:\mathbf{C_2}({\rm proj} \,H) \rightarrow \mathbf{C_n}({\rm proj} \,H)$  defined as follows: for $X=X^1 \rightarrow X^2 \in \mathbf{C_2}({\rm proj} \,H)$ we have $\varphi_i(X)=X[-i]$ where $(X[-i])^j=0$ if $j \notin \{i+1,i+2\}$, $(X[-i])^{i+1}=X^1$ and $(X[-i])^{i+2}=X^2$;  for a morphism $f:X \rightarrow Y \in \mathbf{C_2}({\rm proj} \,H)$ we defined $\varphi_i(f):X[-i]\rightarrow Y[-i] \in \mathbf{C_n}({\rm proj} \,H)$ as the morphism induced by $f$. If $\mathbf{C_2}({\rm proj} \,H )[-i]$ is the image of $\varphi_i$ then we have a bijective correspondence between the indecomposable complexes  of $\mathbf{C_n}({\rm proj} \,H)$ and
the complexes of
$\cup_{i=0}^{n-2}\, \mbox{ind}\, {\mathbf{C_2}({\rm proj} \,H )[-i]}$ which is not a disjoint union. In fact, in \cite[Lemma 5.2]{H:88} it is  showed that an indecomposable complex in the bounded derived category  of modules, $\mathbf{D}^b(H)$,  is a stalk complex with indecomposable stalk. Since $\mathbf{D}^b(H)$ is equivalent to $\mathbf{K}^b({\rm proj} \,H)$ we get that
an indecomposable complex  $X$ in   $\mathbf{K}^b({\rm proj} \,H)$ is of the form
$X:0 \rightarrow \overset{j-1}{\cdots}\rightarrow  X^{j} \stackrel{d}{\rightarrow} X^{j+1}  \rightarrow \overset{n-j-1}{\cdots}\rightarrow 0 $  with  $d$ a monomorphism and $ X^{j+1} \neq 0$ or either of the form $S(P) : P \rightarrow 0\rightarrow \dots \rightarrow  0 \rightarrow 0 $.  By \cite[Corollary 5.7]{BSZ:04} we know that $\mathfrak{L}_n \simeq \overline{\mathbf{C_n}}({\rm proj} \,H)$. Then, the indecomposable complexes in $\mathbf{C_n}({\rm proj} \,H)$ which are not injective objects are the same that the objects in $\mathfrak{L}_n$. Moreover, by \cite[Proposition 6.6]{BSZ:04} the differential $d$ satisfies that $d(X^j) \subset \mbox{rad}\;X^{j+1}$ where $\mbox{rad}$ is the Jacobson radical of $H$. The injective complexes in $\mathbf{C_n}({\rm proj} \,H)$ were defined above.

Moreover, by \cite[Proposition 3.3]{B:04} we know that ${\rm mod}\,H \simeq \overline{\mathbf{C_2}}({\rm proj} \,H)$ and by  \cite{H:88} we know the correspondence between the irreducible morphisms between indecomposable objects in $\mathbf{K}^b({\rm proj} \,H)$ and ${\rm mod}\,H$.
Hence, the embedding $\varphi_i$ satisfies that if $f:X \rightarrow Y \in \mathbf{C_2}({\rm proj} \,H)$ is irreducible then $\varphi_i(f):X[-i]\rightarrow Y[-i]$ is irreducible in $\mathbf{C_n}({\rm proj} \,H)$. Then, the Auslander-Reiten quiver of $\mathbf{C_n}({\rm proj} \,H)$ is equal to $n-1$ copies of the Auslander-Reiten quiver of $\mathbf{C_2}({\rm proj} \,H)$. Moreover, $\mbox{ind}\, {\mathbf{C_2}({\rm proj} \,H )[-j]} \cap \mbox{ind}\, {\mathbf{C_2}({\rm proj} \,H )[-(j+1)]} \neq \emptyset$ for $j=0, \dots, n-3$ with $n\geq 3$, and we conclude that the quiver is connected.
\end{sinob}

\begin{sinob}

We say that a morphism $f$ is a source morphism or a  minimal left almost split morphism if the following conditions hold:
$(i)$ $f$ is not a section, $(ii)$  for each morphism $g:X \rightarrow M$ that is not a section $g$ factors through $f$ and,
$(iii)$ for every morphism $h:Y \rightarrow Y$, if $f=hf$ then $h \in \mbox{Aut}(Y)$.

Dually, we define a sink morphism or a  minimal right almost split morphism.

A  sequence $X \stackrel{f} {\rightarrow} Y \stackrel{g} {\rightarrow} Z$ in $\mathbf{C_n}({\rm proj}\, \Lambda)$ that is exact  in $\mathbf{C_n}({\rm mod}\, \Lambda)$ is called an almost split sequence in $\mathbf{C_n}({\rm proj}\, \Lambda)$ if $f$ and $g$ are a source and a sink morphism, respectively.
The object $X$ is called the translate of $Z$ and it is denoted by $A_n(Z) =X$.
In \cite{BSZ:04}, the authors characterized the translate of $Z$ and proved the existence of almost split sequences  in $\mathbf{C_n}({\rm proj}\, \Lambda)$.
\end{sinob}

\begin{sinob}
By $\Re_{\mathbf{C_{n}}({\rm proj}\, \Lambda)}$ we denote the radical of the category ${\mathbf{C_{n}}({\rm proj}\, \Lambda)}$, which is the ideal generated  by the non-isomorphisms between indecomposable objects.
For a natural number $m$ we denote the powers of the radical by  $\Re^m_{\mathbf{C_{n}}({\rm proj}\, \Lambda)}$. Finally, by the infinite radical $\Re^{\infty}_{\mathbf{C_{n}}({\rm proj}\, \Lambda)}$ we mean the intersection of all the powers $\Re^i_{\mathbf{C_{n}}({\rm proj}\, \Lambda)}$ for  $i \geq 1$.

In \cite{CPS1}, the authors of this paper determined how to knit the Auslander-Reiten quiver of ${\mathbf{C_n}({\rm  proj}\, \Lambda)}$, $\Gamma_{\mathbf{C_n}({\rm  proj}\, \Lambda)}$.

Following \cite{SK}, we say that a connected component  $\Gamma$ of $\Gamma_{\mathbf{C_{n}}({\rm proj}\, \Lambda)}$ is generalized standard if  $\Re^{\infty}_{\mathbf{C_{n}}({\rm proj}\, \Lambda)}(X,Y)=0$ for all $X$ and $Y$ in $\Gamma$.
\end{sinob}

\begin{sinob}\label{GM}
By  \cite[Corollary 2, Proposition 3]{GM:08} if
  $f = \{f^i\}_{i=1}^n:X \rightarrow Y$
is an irreducible
morphism in $\mathbf{C_n}({\rm proj} \,\Lambda)$ then, one of the
following statements hold.

 $(sec)$ For each $i \in \{1,...,n\}$, the morphisms $f^i$ are sections in ${\rm proj} \,\Lambda.$

 $(ret)$ For each $i \in \{1,...,n\}$, the morphisms $f^i$ are retractions in ${\rm proj} \,\Lambda.$

 $(ret-irred-sec)$ There exists $i \in \{1,...,n\}$ such that $f^i$ is irreducible in ${\rm proj}\,\Lambda$, the morphisms $f^j$ are sections for all $j>i$ and the morphisms $f^j$ are retractions for $j<i.$
\end{sinob}

\begin{sinob}\label{T(S)} Let   $P$  be  an indecomposable projective $\Lambda$-module and  $\cdots \overset{d^{-3}}{\rightarrow}
R^{-2}\overset{d^{-2}}{\rightarrow} R^{-1}
\overset{d^{-1}}{\rightarrow} P \rightarrow P/{\rm rad}P
\rightarrow 0$ be a minimal projective resolution of  $P/{\rm
rad}P$ in ${\rm mod}\; \Lambda$ and $P/{\rm rad}P \rightarrow
I^0 \overset{g^0}{\rightarrow} I^1 \overset{g^1}{\rightarrow}
\cdots \overset{g^{n-j-1}}{\rightarrow} I^{n-j} \rightarrow
\cdots$ a minimal injective co-resolution  of $P/{\rm
rad}P$.  Let $R^0=P$ and consider $R=(R^j,d^j_R)$ with $R^j=0$ if $j> 0$ and $d^j_R=d^j$ for $j < 0$. We define $R_j(P)= F(R[-j-1])$. Similarly. let  $L^0=P$ and consider $L=(L^s,d_L^s)$, with $L^s = \nu^{-1}(I^s)$, for  $s \geq 1$, where $\nu$ is the  Nakayama functor, $L^s=0$ for $j< 0$ and $d_L^s=\nu^{-1}(g^s).$  We define $L_j(P)=G(L[-j])$.

Now, consider the following  morphisms of complexes in $\mathbf{C_n}({\rm proj} \,\Lambda)$.

\begin{displaymath}
\def\objectstyle{\scriptstyle}
\def\labelstyle{\scriptstyle}
   \xymatrix   @R=3.5mm { R_j(P)) \ar[d]^{\sigma_j^n} &  & R^{-j}\ar[d]^{id}
   \ar[rr]^{d^{-j}}
 &   &  \cdots \,\,\,  R^{-1}\ar[rr]^{d^{-1}} \ar[d]^{id} & & P \ar[rr]\ar[d]^{d^0_L} & & 0 \ar[rr]\ar[d]& & \cdots \,  0 \,\,\,\, \ar[d]\\
  B_j(P) \ar[d]^{\tau_j^n} & & R^{-j}
   \ar[rr]^{d^{-j}} \ar[d] & & \cdots \,\,\, R^{-1} \ar[rr]^{d^0_L d^{-1}} \ar[d]^{d^{-1}} & & L^1 \ar[rr]^{d^1_L}\ar[d]^{id} &  & L^2 \ar[rr]^{d^2_L} \ar[d]^{id} & & \cdots \,\,\, L^{n-j} \ar[d]^{id} \\
 L_j(P)  & & 0
   \ar[rr] & & \cdots \,\,\,  P\,\,\,\, \ar[rr]^{d^0_L} & & L^1 \ar[rr]^{d^1_L} &  & L^2 \ar[rr]^{d^2_L} & & \cdots \,\,\, L^{n-j} }
\end{displaymath} \noindent and

\begin{displaymath}
\def\objectstyle{\scriptstyle}
\def\labelstyle{\scriptstyle}
   \xymatrix  @R=3.5mm  { R_j(P) \ar[d]^{\rho_j^n} &  & R^{-j}\ar[d]^{}
   \ar[rr]^{d^{-j}}
 &   &  \cdots \,\,\,  R^{-1}\ar[rr]^{d^{-1}} \ar[d]^{d^{-1}} & & P \ar[rr]\ar[d]^{id} & & 0 \ar[rr]\ar[d]& & \cdots \,  0 \,\,\,\, \ar[d]\\
  J_j(P) \ar[d]^{\lambda_j^n} & & 0
   \ar[rr]^{} \ar[d] & & \cdots \,\,\, P \ar[rr]^{id} \ar[d]^{id} & & P \ar[rr]^{}\ar[d]^{d_L^0} &  & 0 \ar[rr]^{} \ar[d]^{} & & \cdots \,\,\, 0 \ar[d]^{} \\
 L_j(P)  & & 0
   \ar[rr] & & \cdots \,\,\,  P\,\,\,\, \ar[rr]^{d^0_L} & & L^1 \ar[rr]^{d^1_L} &  & L^2 \ar[rr]^{d^2_L} & & \cdots \,\,\, L^{n-j} }
\end{displaymath}

\vspace{.1in}

\noindent In  \cite[Proposition 8.5, 8.7, and 8.8]{BSZ:04} the authors proved that
$\lambda _{j}^n: J_{j}(P)\rightarrow L_{j}(P)$
\noindent is a minimal left almost split morphism and
$\rho_j^n:R_j(P) \rightarrow J_j(P) $
\noindent  is a minimal right almost split morphism  for $j=1, \cdots, n-1$. Moreover, they showed that
$R_{j}(P) \stackrel{(\rho _{j}^n,\sigma _{j}^n)^{t}} \rightarrow J_{j}(P) \oplus
B_{j}(P) \stackrel{(\lambda _{j}^n,\tau _{j}^n)} \rightarrow L_{j}(P)$  is an almost split sequence in $\mathbf{C_n}(\rm{proj}
\,\Lambda)$. Therefore, the morphisms $\sigma^n_j:R_j(P) \rightarrow B_j(P)$ and $\tau_j^n:B_j(P) \rightarrow L_j(P)$ are irreducible. The complexes
$R_j(P)$, $J_j(P)$ and  $L_j(P)$ are alway indecomposable but the complex $B_j(P)$ may not be indecomposable.
\end{sinob}

\begin{sinob}
If $f= \{f^i\}_{i=1}^n: X \rightarrow Y$  is a morphism in
$\mathbf{C_n}(\rm{mod} \,\Lambda)$ then ${\rm Ker} \,f$ (${\rm Coker} \,f $) is the complex $({\rm Ker}\,
f^i, d_{{\rm Ker}\,f^i})$ (${\rm Coker} \,f =
({\rm Coker}\,f^i, d'_i)$), where $d_{{\rm Ker}\,f^i}$ is the restriction of $d_{X^i}$ to ${\rm Ker}\,f^i$ ($d'_i$ is induced by the morphisms
$f^i$ and $f^{i+1}$ from  ${\rm Coker}\,f^i$ to ${\rm
Coker}\,f^{i+1}$, respectively).

By Lemma \ref{Coker}, we know that the kernels and cokernels in $\mathbf{C_n}({\rm proj} \,\Lambda)$ do not belong necessarily to the category, since ${\rm proj} \,\Lambda$ is not closed under kernels and cokernels. Hence, in general $\mathbf{C_n}({\rm proj} \,\Lambda)$  is not exact abelian.
\end{sinob}

\section{On the Kernel and Cokernel of irreducible morphisms in
$\mathbf{C_n}(\rm{proj}\, \Lambda)$}

In a module category, the kernel (cokernel) has played an important role in order to determine if an irreducible morphism is of finite left (right) degree. For such a reason, we are interested to study the kernel (cokernel) of irreducible morphisms in $\mathbf{C_n}(\rm{proj} \,\Lambda)$, for $\Lambda$ an artin algebra. Moreover, we are also interested in knowing when the kernel (cokernel) of an irreducible morphism belongs to $\mathbf{C_n}(\rm{proj} \,\Lambda)$, since such a category  is not abelian.
\vspace{.05in}

We start proving some necessary lemmas.

\begin{lem}\label{Coker} Let $\Lambda$ be an artin algebra and $f:X \rightarrow Y$ be
an irreducible morphism in $\rm{proj}\,\Lambda$. Then, ${\rm Coker}\,f
\notin \rm{proj}\,\Lambda$.
\end{lem}

\begin{proof} Given $f:X \rightarrow Y$
irreducible in $\rm{proj}\,\Lambda$ then ${\rm Coker}\,f \neq 0$, otherwise $f$ splits. Assume that ${\rm
Coker}\,f \in \rm{proj}\,\Lambda$. Then, $Y \simeq {\rm
Coker}\,f \oplus {\rm Im}f$. Hence, ${\rm Im}f \in \rm{proj}\,\Lambda$.
Consider $g:X \rightarrow {\rm Im}f $ and  $h: {\rm
Im}f \rightarrow Y$ such that $g(x)=f(x)$ for all $x \in X$ and where $h$ is the inclusion morphism. Hence, $f= h g$. Then, either $g$ is a section or $h$ is a retraction. In the former case there exists $g': {\rm Im}f \rightarrow X$ such that $g'g=id_X$. If we consider the morphism $(0,g'):{\rm
Coker}\,f \oplus {\rm Im}f \rightarrow X$ then $(0,g')f=id_X$, getting a contradiction.
In the latter case, we infer that $f$ is an epimorphism, a contradiction. Therefore, ${\rm Coker}\,f \notin \rm{proj}\,\Lambda$.
\end{proof}

We recall that a morphism in ${\rm proj}\,H,$ is irreducible in ${\rm mod}\,H$ if and only if it is irreducible in ${\rm proj}\,H,$ whenever $H$ is a hereditary algebra, see \cite[Lemma 3.3]{CPS1}.

As a consequence of Lemma \ref{Coker} and the above comment we get the following result.

\begin{coro}\label{Cokernoproy} Let $H$ be a hereditary artin algebra and $f:X \rightarrow Y$ be an irreducible morphism in ${\rm proj}\,{H}$. Then, the following conditions hold.
\begin{enumerate}
\item[(a)] The morphism $f$ is a monomorphism.
\item[(b)]  There is not a non-zero morphism $g:Y \rightarrow P$ in
${\rm proj}\,{H}$ such that $g f=0$.
\end{enumerate}
\end{coro}
\begin{proof} $(a).$ Follows from the fact that if $f:X \rightarrow Y$ is irreducible in ${\rm proj}\,{H}$ then $f$ is irreducible in ${\rm mod}\,{H}$. By Lemma \ref{Coker} $f$ is not an epimorphism. Therefore, $f$ is a monomorphism in ${\rm mod}\,{H}$.

$(b).$ Assume that there exists a non-zero morphism $g:Y \rightarrow P$
in ${\rm proj}\,{ H}$ such that $g f=0$. Then, there exists $g':Y \rightarrow {\rm Coker}\,f$ and  $h: {\rm
Coker}\,f \rightarrow P$ in ${\rm mod}\,{H}$ such that $g =h g'$. Since
$H$ is hereditary then ${\rm Im}\, h \in {\rm proj}\,{H}$ and therefore ${\rm Im}\, h$ is a direct summand of ${\rm Coker}\,f$.
Moreover, by  \cite[V,  Proposition 5.6]{ARS:95} since ${\rm Coker}\,f $ is indecomposable, we have that ${\rm Im}\, h \simeq{\rm Coker}\,f \in  {\rm proj}\,{H}$, contradicting Lemma \ref{Coker}.
\end{proof}

\begin{rem}\label{Ker 1} {\rm Let $f:X \rightarrow Y$ be an irreducible morphism  in
$\mathbf{C_n}(\rm{proj} \,\Lambda)$ with $X$ or $Y$
indecomposable. If all entries of  $f$ are retractions then ${\rm Ker} \,f \neq 0$, otherwise $f$ is an isomorphism. Moreover, ${\rm Ker} \,f \in \mathbf{C_n}(\rm{proj}
\,\Lambda)$. In fact, since for each  $i$
the exact sequence $0 \rightarrow {\rm Ker} \,f^i \rightarrow X^i
\rightarrow Y^i \rightarrow 0$ splits then ${\rm Ker} \,f^i$
is a direct summand of $X^i$ and, therefore ${\rm Ker} \,f^i$ is
projective.  The converse is not true. That is, if ${\rm Ker} \,f \neq 0$ then the morphism $f$ does not  has necessarily retractions in all its entries.}
\end{rem}

It is well known that an irreducible morphism in
$\mbox{mod}\,\Lambda$ is either a monomorphism or an epimorphism.
In general, this property does not hold in $\mathbf{C_n}({\rm proj} \,\Lambda)$. In case
$\Lambda$ is hereditary then we shall prove that the irreducible morphisms have such a property.

\begin{prop}\label{moe}
Let $\Lambda$ be an artin algebra and $f=\{f^i\}_{i=1}^n: X \rightarrow
Y$ be an irreducible morphism in $\mathbf{C_n}({\rm proj}
\,\Lambda).$  If $f$ is not an epimorphism and ${\rm Im}f^i \in
{\rm proj} \,\Lambda$ for $1 \leq i \leq n$ then $f$ is a
monomorphism.
If $H$ is a hereditary algebra then $f$ is either a
monomorphism or an epimorphism.
\end{prop}

\begin{proof}
Consider ${\rm Im}f =({\rm Im} f^i, d^i_{{\rm Im} f})$ where $d^i_{{\rm Im} f}$ is the restriction morphism  $d^i_{Y}$
to ${\rm Im}f^i$ for $i=1, \dots, n$. Assume that $f$ is not an epimorphism and that ${\rm
Im}f \in \mathbf{C_n}({\rm proj} \,\Lambda)$. Then, $f$ can be written as the composition $\beta
\alpha$ where $\alpha: X \rightarrow {{\rm Im}f}$ is the co-restriction morphism and $\beta: {\rm
Im}f {\rightarrow} Y$ is the inclusion morphism. Furthermore, $\alpha$  ia an epimorphism. Hence, since $f$ is not an
epimorphism, then $\beta$ is not an epimorphism. Moreover, $\beta$
does not a  split epimorphism. Hence, since $f$ is irreducible then $\alpha$ is a split
monomorphism. Therefore, $f$ is a monomorphism.

Now, if $H$ is a hereditary algebra then for each $i=1, \dots, n$  we have that ${\rm Im} f^i \in {\rm proj}\, H$. Then, $f$ is either
an epimorphism or a monomorphism.
\end{proof}

\begin{prop}\label{Ker en Cn(H)} Let $H$ be a hereditary artin algebra and $f=\{f^i\}_{i=1}^n:X \rightarrow Y$ be an irreducible morphism in
$\mathbf{C_n}({\rm proj}\,{ H})$ with $X$ or $Y$ indecomposable. Then,
\begin{enumerate}
\item [(1)] The morphism  $f$ is of type $(ret)$ if and only if ${\rm Ker}
\,f \in \mathbf{C_n}({\rm proj}\,{ H}) \backslash \{0\}.$
\item [(2)] If  $f$ is of type $(ret-irred-sec)$  then ${\rm Ker}\,f =0.$
\end{enumerate}
\end{prop}

\begin{proof} $(1)$. If $f$ is of type $(ret)$, then ${\rm Ker}f \neq 0$  by Remark \ref{Ker 1}. Conversely, assume that
${\rm Ker} \,f \in \mathbf{C_n}({\rm proj}\,{ H})\backslash \{0\}$. By Proposition \ref{moe} we have that $f$ is an epimorphism. Then, each entry $f^i$ is a retraction. Therefore, $f$ is of type $(ret)$.

$(2)$ If $f$ is of type $(ret-irred-sec)$ then by Proposition \ref{moe} we claim that  $f$ is a
monomorphism. Otherwise, $f$ is an epimorphism and hence, the component of $f$ that is irreducible splits. Therefore, ${\rm Ker}\,f = 0$.
\end{proof}

In our next result we show that the cokernel of an irreducible morphism determines the shape of the morphism. More precisely, we prove the following result.

\begin{prop}\label{Nu-Conucleo} Let $\Lambda$ be an artin algebra and  $f:X \rightarrow Y$ be an irreducible morphism in
$\mathbf{C_n}(\rm{proj} \,\Lambda)$ with  $X$ or $Y$ indecomposable. Then, the following conditions hold.

\begin{enumerate}
\item [(a)] The morphism $f$ is of type $(sec)$ if and only if ${\rm Coker}
\,f \in \mathbf{C_n}(\rm{proj} \,\Lambda)\backslash \{0\}$.
\item [(b)]  The morphism $f$ is of type  $(ret-irred-sec)$ if and only if  ${\rm Coker} \,f \notin \mathbf{C_n}(\rm{proj} \,\Lambda)$ \rm{(in particular, $\mathbf{C_n}(\rm{proj} \,\Lambda)\neq 0$)}.
\item [(c)] The morphism $f$ is of type  $(ret)$ if and only if  ${\,\rm Coker}\,f =0.$
\end{enumerate}
\end{prop}

\begin{proof}
$(a).$ If all entries of  $f$ are sections then ${\rm Coker}\,f^i \in {\rm proj}\,{\Lambda}$, since for each $i$ the exact sequence $0
\rightarrow X^i \rightarrow Y^i \rightarrow {\rm Coker} \,f^i
\rightarrow 0$ splits. Conversely, if $0 \neq {\rm Coker}
\,f \in \mathbf{C_n}(\rm{proj} \,\Lambda)$ then $f$ is not of type $(ret)$. Assume that $f$ is of type $(ret-irred-sec)$ with $f^i:X^i \rightarrow Y^i$ irreducible in $\rm{proj} \,\Lambda$. By Lemma \ref{Coker} we have that ${\rm Coker} \,f^i \notin \rm{proj}
\,\Lambda$. Then, we infer by (\ref{GM}) that $f$ is of type $(sec)$.

$(b).$ Let $f$ be a morphism with $f^j$ irreducible
in $\rm{proj} \,\Lambda$. By Lemma \ref{Coker}, ${\rm Coker} \,f^j \notin \rm{proj}
\,\Lambda$. Conversely, if ${\rm Coker} \,f \notin \mathbf{C_n}(\rm{proj} \,\Lambda)$ then $f$ is not of type $(ret)$. By $(a)$, $f$ is of type $(ret-irred-sec)$.

$(c).$ Follows trivially.
 \end{proof}

Next we show an example of an  irreducible morphism $f$ in $\mathbf{C_2}(\rm{proj} \,\Lambda)$ such that ${\rm Ker}\,f$ belongs to $\mathbf{C_2}(\rm{proj} \,\Lambda)$ and it is not indecomposable.

\begin{ej}\label{ker descomponible} {\em Consider the path algebra given by the quiver}

$$
\xymatrix{ 3 & 2 \ar[l]_{\alpha} & 1 \ar[r]^{\gamma}
\ar[l]_{\beta} & 4 \ar[r]^{\delta} & 5}$$

\noindent {\em with
$\beta\alpha=0= \delta\gamma$.  The Auslander-Reiten quiver of
$\mathbf{C_2}(\rm{proj} \;\Lambda)$ is the following}

 \tiny
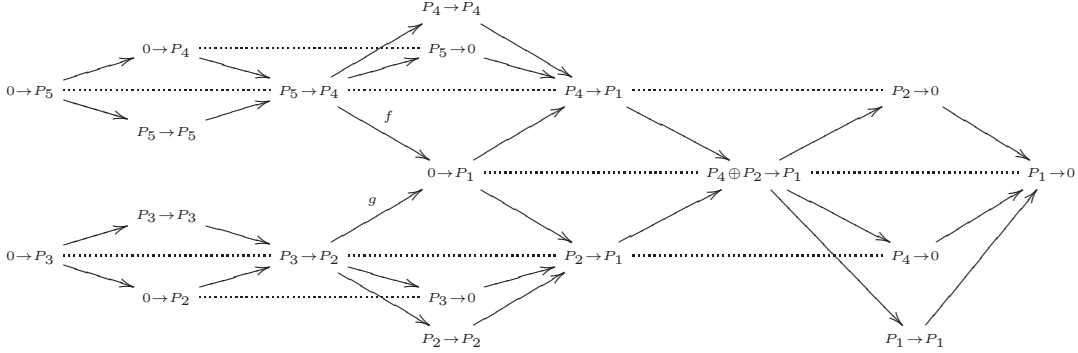
\begin{displaymath}
\def\objectstyle{\scriptstyle}
\def\labelstyle{\scriptstyle}
   \xymatrix @R=1.5mm  {& & & P_4 \rightarrow P_4 \ar[rdd] & & & &  \\
   & 0\rightarrow P_4\ar@{.}[rr]\ar[rd] & & P_5 \rightarrow 0 \ar[rd]&  & &  & \\
0\rightarrow P_5 \ar[ru]\ar[rd]\ar@{.}[rr] & & P_5\rightarrow P_4 \ar[rdd]^{ f}\ar[ru]\ar[ruu]\ar@{.}[rr] & &
P_4\rightarrow P_1 \ar[rdd]\ar@{.}[rr]& & P_2\rightarrow 0 \ar[rdd]&  \\ & P_5\rightarrow P_5
\ar[ru]& & & & & & \\
 & & & 0\rightarrow P_1 \ar[rdd]\ar[ruu]\ar@{.}[rr]& & P_4\oplus P_2\rightarrow P_1 \ar[ruu]\ar[rdddd]\ar[rdd]\ar@{.}[rr]& & P_1\rightarrow 0 \\
& P_3\rightarrow P_3\ar[rd]&  & & & &  & \\
0\rightarrow P_3\ar[rd]\ar[ru]\ar@{.}[rr]& & P_3\rightarrow P_2\ar[rd]\ar[rdd]\ar[ruu]^{ g} \ar@{.}[rr]& &P_2\rightarrow P_1 \ar[ruu]\ar@{.}[rr]& & P_4\rightarrow 0 \ar[ruu]& \\
& 0\rightarrow P_2\ar[ru]\ar@{.}[rr]&  & P_3\rightarrow 0\ar[ru]& & &  & \\
& &  & P_2\rightarrow P_2\ar[ruu]& & &  P_1\rightarrow P_1\ar[ruuuu]& }
\end{displaymath}

\normalsize {\em We observe that the kernel of the irreducible morphism  $(f , g): (P_5 - P_4) \oplus
(P_3-P_2)\rightarrow (0-P_1)$ is the complex $(P_5-P_5) \oplus (P_3-P_3)$ which is not indecomposable.}
\end{ej}

Next, we prove that all irreducible morphisms in $\mathbf{C_n}(\rm{proj} \,\Lambda)$ of type $(sec)$ ($(ret)$, respectively) have indecomposable cokernel (kernel, respectively) in $\mathbf{C_n}(\rm{proj} \,\Lambda)$.

\begin{thm} \label{Nu yconu indesc} Let $\Lambda$ be an artin algebra and $f:X \rightarrow Y$ be an irreducible morphism  in
$\mathbf{C_n}(\rm{proj} \,\Lambda)$ with  $X$ or  $Y$ indecomposable. Then, the following conditions hold.
\begin{enumerate}
\item [(a)] If $f$ is of type $(ret)$ then ${\rm Ker}\,f$ is  indecomposable and ${\rm Ker}\,f \in \mathbf{C_n}(\rm{proj} \,\Lambda)$.
\item [(b)] If $f$ is of type $(sec)$ then ${\rm Coker}\,f$ is  indecomposable and ${\rm Coker}\,f \in \mathbf{C_n}(\rm{proj} \,\Lambda)$.
\end{enumerate}
\end{thm}

\begin{proof} We only prove $(a)$, since $(b)$ follows similarly.

$(a).$ Let $f=\{f^i\}_{i=1}^{n}$ be of type $(ret)$. By Remark \ref{Ker 1} we know that ${\rm Ker}\,f$ is not zero and belongs to $\mathbf{C_n}(\rm{proj} \,\Lambda)$. Consider $K= {\rm Ker}\,f$ and  $g=\{g^i\}_{i=1}^{n}: K \rightarrow X$ the inclusion morphism. Assume that $K = N \oplus T$ with $N$ and $T$ non-zero complexes. Let  $\pi_{_N}:K \rightarrow N$ be the projection morphism. Since for each $i$, $g^i$ is a section and $f^i$ is a retraction then without loss of generality we may write $g^i= (1,0)^t$,  $f^i= (0,1)$ and $X^i= K^i \oplus Y^i$. Hence, we have the following commutative diagram:

\begin{displaymath}
\xymatrix   {K \ar[d]^g & &\ar[r] & K^i \ar[r]^{d_K^i} \ar[d]^{(1,0)^t}    & K^{i+1} \ar[r] \ar[d]^{(1,0)^t} &  \\ X \ar[d]^f & & \ar[r]  &  K^i \oplus Y^i \ar[r]^{d_X^i} \ar[d]^{(0,1)} & K^{i+1} \oplus Y^{i+1} \ar[d]^{(0,1)}\ar[r] &
             \\ Y & & \ar[r] & Y^i \ar[r]^{d_Y^i} & Y^{i+1} \ar[r] &
              }
\end{displaymath}
\vspace{.1in}

\noindent with  $d_X^i= \left( \begin{array}{cc} d_K^i & b^i \\  0 & d_Y^i \\ \end{array}\right)$, $b^i: Y^i \rightarrow K^{i+1}$ in
${\rm proj}\,\Lambda$ and where $g:K \rightarrow X$ is the inclusion morphism. Since $d_X^{i+1} d_X^i =0$ then $d_K^{i+1}b^i + b^{i+1}d^i_Y=0.$

Now, consider $d: K \oplus Y \rightarrow E$, where $d^i =\left( \begin{array}{cc} \pi_N^i & 0 \\  0 & 1 \\ \end{array}\right)$,
$E=(E^i,d_E^i)$,
$E^i =N^i \oplus Y^i$ and
 $d_E^i= \left( \begin{array}{cc} d_N^i & \pi_N^{i+1}b^i \\  0 & d_Y^i \\ \end{array}\right)$. Note that $E$ is a complex. In fact, $d^{i+1}_E d^i_E = \left( \begin{array}{cc} 0 & d_N^{i+1}\pi_N^{i+1}b^i + \pi_N^{i+2}b^{i+1}d^i_Y \\  0 & 0 \\ \end{array}\right)$ and, since $d_K^{i+1}b^i + b^{i+1}d^i_Y=0$ then

$$\begin{array}{lcl}
 d_N^{i+1}\pi_N^{i+1}b^i + \pi_N^{i+2}b^{i+1}d^i_Y
& = d_N^{i+1}\pi_N^{i+1}b^i - \pi_N^{i+2}d_K^{i+1}b^i \\
 & = (d_N^{i+1}\pi_N^{i+1} - \pi_N^{i+2}d_K^{i+1})b^i \\  & \hspace{-1.4in}= 0.
\end{array}$$

\noindent Moreover, for each $i$ we have the following commutative diagram:

\begin{displaymath}
\def\objectstyle{\scriptstyle}
\def\labelstyle{\scriptstyle}
\xymatrix @!0 @C=9mm @R=9mm  { &   & K^{i+1}\ar[rrr]^{(1,0)^t}\ar[ddd]^{\pi_{_N}^{i+1}}  & & &  K^{i+1} \oplus Y^{i+1}\ar[rrr] ^{(0,1)} \ar[ddd]^{d^{i+1}} & & & Y^{i+1} \ar@{=}[ddd] \\ & & & & & & & &  \\
  K^i \ar[rrr]^{(1,0)^t} \ar[ddd]_{\pi_N^i}\ar[rruu]^{d_K^i}& & & K^i \oplus Y^i  \ar[rrr]^{(0,1)} \ar[ddd]_{d^i} \ar[rruu]^{d_X^i}& & & Y^i \ar@{=}[ddd] \ar[rruu]^{d_Y^i}& & \\ & & N^{i+1} \ar[rrr]^{(1,0)^t}  & & & E^{i+1}\ar[rrr]^{(0,1)}  & & & Y^{i+1} \\  & & & & & & & &
  \\ N^i \ar[rrr]^{(1,0)^t}\ar[rruu]^{d_N^i} &  & &  E^i \ar[rrr]^{(0,1)}\ar[rruu]^{d_E^i} & & & Y^i \ar[rruu]^{d_Y^i} & &}
\end{displaymath}
\vspace{.1in}

\noindent Let  $g'= (1, 0)^t:N \rightarrow E$ and $f'=(0,1):E\rightarrow Y$. Then, the following diagram
$$
\xymatrix{K \ar[d]^{\pi_N}\ar[r]^{g} & X \ar[d]^{d} \ar[r]^{f}  & Y \ar@{=}[d]\\
 N \ar[r]^{g'} & E \ar[r]^{f'} & Y} $$ commutes.
Since $f=f'd$ then $d$ is a section or else $f'$ a retraction. In the former case there exists $t:E \rightarrow X$ such that $t d= id_{_X}$. Then, $t g' \pi_{_N} = t d g = id_{_X} g = g$ and, because $g$ is a monomorphism then $\pi_{_N}$ is an isomorphism, contradicting the assumption that $N$ and $T$ are non-zero. In the latter case, we infer that $g'$ is a section and, there exists $h_{_N}: E \rightarrow N$ such that $h_{_N} g'=id_{_N}$. Hence, $h_{_N} d g = h_{_N} g' \pi_{_N} = id_{_N} \pi_{_N}= \pi_{_N}$.

Analogously, let $\pi_{_T}:K \rightarrow T$  be the projection morphism and consider $d':K \oplus Y \rightarrow E'$, where $d'^i =\left( \begin{array}{cc} \pi_T^i & 0 \\  0 & 1 \\ \end{array}\right)$,
$E'=(E'^i,d_E^i)$,
$E'^i =T^i \oplus Y^i$ and
 ${d'}_E^i= \left( \begin{array}{cc} d_T^i & \pi_T^{i+1}b^i \\  0 & d_Y^i \\ \end{array}\right)$. Without loss of generality, we may write $\widetilde{g} =(1,0)^t:T \rightarrow E'$ and $\widetilde{f}
=(0,1):E' \rightarrow Y$ with $f=\widetilde{f} d'$. As we analyzed before, we get that  $d'$ is not a section. Then $\widetilde{f}$ is a retraction and,  there is $h_{_T}:E' \rightarrow T$ such that $h_{_T}\widetilde{g}=id_{_T}$ and, $h_{_T}d' g= \pi_{_T}$.
Hence, $\left( \begin{array}{c} h_{_N} d \\  h_{_T} d' \\ \end{array}\right) g =id_{_K}$ and, therefore $g$ is a section. Then, $f$ is a retraction, a contradiction. Thus, ${\rm Ker}\,f$ is indecomposable.
\end{proof}

\section{Main results}
It is well-known that an artin algebra  $\Lambda$ is representation-finite if and only if the infinite radical of $\rm{mod}\,\Lambda$ vanishes, see  \cite[V, Sec. 7]{ARS:95}. Using similar arguments as the ones to prove the mentioned result, we affirm that the result also holds for the category $\mathbf{C_n}(\rm{proj}\,\Lambda)$. In fact, if $M$ is a direct sum of copies of a complete set of representatives of isomorphism classes of indecomposable objects in $\mathbf{C_n}(\rm{proj}\,\Lambda)$, then any radical morphism in $\mathbf{C_n}(\rm{proj}\,\Lambda)$ between indecomposable complexes lies in the radical of the artin algebra $\mbox{End}(M)$. Accordingly, the radical of $\mathbf{C_n}(\rm{proj}\,\Lambda)$ is nilpotent, and hence the infinite radical is zero.
\vspace{.1in}

Let $f:X \rightarrow Y$ be an irreducible morphism in $\mathbf{C_n}(\rm{proj} \,\Lambda)$ with $X$ or $Y$
indecomposable. Following \cite{L:10}, we define the left degree of $f$, $d_l(f)$, to be the least integer $m \geq 1$, such that  there is a morphism $\varphi \in \Re^m_{\mathbf{C_n}(\rm{proj} \,\Lambda)}(M,X) \backslash
\Re^{m+1}_{\mathbf{C_n}(\rm{proj} \,\Lambda)}(M,X)$ such that $f \varphi \in \Re^
{m+2}_{\mathbf{C_n}(\rm{proj} \,\Lambda)}(M,Y)$ for some $M \in \mathbf{C_n}(\rm{proj} \,\Lambda)$; and it is infinite if no such integer exist.
Dually, we define the right degree of $f$, $d_r(f)$.
\vspace{.05in}

We recall that if $X$ and $Y$ are indecomposable complexes in $\mathbf{C_n}({\rm proj}\, \Lambda)$ then ${\rm Irr}_{\mathbf{C_n}({\rm proj}\, \Lambda)}(X,Y)= \Re_{\mathbf{C_n}({\rm proj}\, \Lambda)} (X,Y) / \Re_{\mathbf{C_n}({\rm proj}\, \Lambda)} ^2 (X,Y)$. We recall that ${\rm Irr}_{\mathbf{C_n}({\rm proj}\, \Lambda)}(X,Y)$ is a $k_Y-k_X-$bimodule, where $k_Z={\rm End}_{\mathbf{C_n}({\rm proj}\, \Lambda)}(Z)/\Re_{\mathbf{C_n}({\rm proj}\, \Lambda)}(Z,Z)$.
\vspace{.05in}

We start by proving the following useful result. The statement and proof follow from \cite[Corollary 1.7]{SL:96}. Moreover, in \cite[Lemma 3.6]{L:10} a similar statement is proved in the setting of left Auslander-Reiten categories.

\begin{prop}\label{val trivial} Let $f:X\rightarrow Y$ and $g:X \rightarrow Y$ be irreducible morphisms between indecomposable complexes in  $\mathbf{C_n}(\rm{proj} \,\Lambda)$, where ${\rm dim}_{k_X}{\rm Irr}(X,Y)={\rm dim}_{k_Y}{\rm Irr}(X,Y)= 1$. Then, $d_l(f)=d_l(g)$
 and $d_r(f)=d_r(g)$.
\end{prop}

\begin{proof} We only prove the result for the left degree. Let $d_l(f)=m$. By definition there exists a complex $Z \in \mathbf{C_n}(\rm{proj} \,\Lambda)$ and  $\varphi \in  \Re^m_{\mathbf{C_n}(\rm{proj} \,\Lambda)}(Z,X) \backslash \Re^{m+1}_{\mathbf{C_n}(\rm{proj} \,\Lambda)}(Z,X)$ such that $f \varphi \in \Re^{m+2}_{\mathbf{C_n}(\rm{proj} \,\Lambda)}(Z,Y)$. Since  ${\rm dim}_{k_X}{\rm Irr}(X,Y)={\rm dim}_{k_Y}{\rm Irr}(X,Y)= 1$ then $g= f\varphi_{_X}  + \mu= \varphi_{_Y} f+ \mu'$, where $\varphi_{_X} \in \rm{Aut}\,( X )$, $\varphi_{_Y}\in \rm{Aut}\,( Y )$ and $\mu, \mu' \in \Re^2_{\mathbf{C_n}(\rm{proj} \,\Lambda)}(X,Y)$. Hence, $g \varphi = \varphi_{_Y} f \varphi + \mu' \varphi$. Then $g \varphi \in \Re^{m+2}_{\mathbf{C_n}(\rm{proj} \,\Lambda)}(Z,Y)$. If $d_l(g)< m$ then $d_l(f)< m$. Therefore, $d_l(g)= m$.
If $d_l(f) = \infty$ the result follows trivially.
\end{proof}

\subsection{Components with length}

Let $\Gamma$ be a connected component of the Auslander-Reiten quiver
$\Gamma_{\mathbf{C_n}(\rm{proj} \,\Lambda)}$.
We recall that  a path in $\Gamma$ is a sequence of irreducible morphisms between indecomposable objects in $\Gamma$ and,
two paths in $\Gamma$ are parallel if they have the same starting point and the same ending point.

\vspace{.1in}

Following \cite{CPT:04}, we say that a component $\Gamma$ of $\Gamma_{\mathbf{C_n}(\rm{proj} \,\Lambda)}$ is with length if parallel paths have the same length and we denote by $\ell(X,Y)$ the length of a path from $X, Y \in \mbox{mod}\;A$. Otherwise, we say that the component is without length. Observe that a  component with length does not have oriented cycles.

The next result proved in \cite[Proposition 3.1]{CPT:04} can be adapted to ${\mathbf{C_n}(\rm{proj} \,\Lambda)}$ with  a similar proof. We only state the result.

\begin{prop}\label{R} Let $\Gamma$ in $\Gamma_{\mathbf{C_n}(\rm{proj} \,\Lambda)}$ be a
generalized standard component with length. Let $X,Y\in
\Gamma$ such that $\ell(X,Y)=n$. Then:

(a) $\Re_{\mathbf{C_n}(\rm{proj} \,\Lambda)}^{n+1}(X,Y)=0$

(b) If $g:X\rightarrow Y$ is a nonzero morphism then $g\in
\Re_{\mathbf{C_n}(\rm{proj} \,\Lambda)}^{n}(X,Y)\setminus \Re_{\mathbf{C_n}(\rm{proj} \,\Lambda)}^{n+1}(X,Y).$

(c) $\Re_{\mathbf{C_n}(\rm{proj} \,\Lambda)}^{j}(X,Y)= \Re_{\mathbf{C_n}(\rm{proj} \,\Lambda)}^{n}(X,Y)$, for each $j=1,...,n$.
\end{prop}

As an immediate consequence of Proposition \ref{R} we get the following result which is fundamental for the next proofs.

\begin{thm}\label{first}
Let $\Lambda$ be an artin algebra  and $\Gamma$ be a generalized standard component of $\Gamma_{\mathbf{C_n}(\rm{proj} \,\Lambda)}$ with length. Let $f:X\rightarrow Y$ be an irreducible morphism with $X, Y \in \Gamma$. Then:
\begin{itemize}
\item [(a)] $d_l(f)=m$ then there is $M \in
\Gamma$ and $g \in \Re_{{\mathbf{C_n}} (\rm{proj} \,\Lambda)}^{m} (M,X) \backslash
\Re_{{\mathbf{C_n}}(\rm{proj} \,\Lambda)}^{m+1}(M,X)$ such that $fg =0$.
\item [(b)] $d_r(f)=m$ then there is $M \in
\Gamma$ and $g \in \Re_{{\mathbf{C_n}} (\rm{proj} \,\Lambda)}^{m} (Y, M) \backslash
\Re_{{\mathbf{C_n}}(\rm{proj} \,\Lambda)}^{m+1}(Y, M)$ such that $gf =0$.
\end{itemize}
\end{thm}

In  Proposition \ref{Nu-Conucleo} and Theorem \ref{Nu yconu indesc} we study when the  kernel or the cokernel of an irreducible morphism in
$\mathbf{C_n}(\rm{proj} \,\Lambda)$ belong to the category. In this case, we get a similar result as
\cite[Theorem C]{CLT}.

\begin{prop}\label{keringama} Let $\Lambda$ be an artin algebra and $\Gamma$ be a generalized standard component
of $\Gamma_{\mathbf{C_n}(\rm{proj} \,\Lambda)}$ with length. Let  $f:X \rightarrow Y$ be an irreducible morphism
in $\mathbf{C_n}(\rm{proj} \,\Lambda)$, with $X,Y \in \Gamma$.
\begin{itemize}
\item [(a)] Assume that ${\rm Ker}\,f \neq 0$ and ${\rm Ker}\,f \in \mathbf{C_n}(\rm{proj} \,\Lambda)$. Then,
$d_l(f) < \infty$ if and only if there exists an indecomposable direct summand of ${\rm Ker}\, f$ that belongs to $\Gamma$.
\item [(b)] If ${\rm Coker}\,f \neq 0$ and ${\rm Coker}\,f \in \mathbf{C_n}(\rm{proj} \,\Lambda)$ then
 $d_r(f) < \infty$ if and only if ${\rm Coker}\,f \in \Gamma$.
\end{itemize}
\end{prop}

\begin{proof} $(a)$. If  $d_l(f) = m$  then there is $Z \in \Gamma$ and
$g \in \Re_{{\mathbf{C_n}} (\rm{proj} \,\Lambda)}^{m} (Z,X) \backslash
\Re_{{\mathbf{C_n}}(\rm{proj} \,\Lambda)}^{m+1}(Z,X)$ such that $f g =0$. Then, $g$ factors through
${\rm Ker}\,f = \oplus_{i =1}^r K_i$, where $K_i$ is indecomposable for $i=1,\cdots, r$. Consider $j=(j_1,
\cdots, j_r): \oplus_{i =1}^r K_i \rightarrow X$ the inclusion morphism. Then, there is a morphism
$t=(t_1, \cdots, t_r)^t: Z \rightarrow {\rm
Ker}\,f $ such that $j t = g$.
Moreover, since $g \in \Re_{{\mathbf{C_n}} (\rm{proj} \,\Lambda)}^{m} (Z,X) \backslash
\Re_{{\mathbf{C_n}}(\rm{proj} \,\Lambda)}^{m+1}(Z,X)$  then $j \not \in \Re_{{\mathbf{C_n}} (\rm{proj} \,\Lambda)}^{\infty} ({\rm Ker} \,f, X)$.
Hence, $j_{_i} \in
\Re_{{\mathbf{C_n}}(\rm{proj} \,\Lambda) }^{s_i}(K_i, X)\backslash \Re_{{\mathbf{C_n}}(\rm{proj} \,\Lambda)}^{{s_i}+1}(K_i, X)$
for some $i\in \{1, \dots, r\}$ where  $s{_i} \leq m$. Hence, $K_i  \in \Gamma$.

Conversely, let ${\rm Ker}\,f= \oplus_{i =1}^r K_i \in
\mathbf{C_n}(\rm{proj} \,\Lambda)$. Since  there is a direct summand of ${\rm Ker}\, f$ that belongs to $\Gamma$
and  $\Gamma$ is a generalized standard component with length, if $j=(j_1, \cdots, j_r): \oplus_{i =1}^r K_i
\rightarrow X$ is the inclusion morphism then there is an integer $i\in \{1, \dots, r\}$ such that $j_{_i} \in
\Re_{{\mathbf{C_n}}(\rm{proj} \,\Lambda)}^{s_i} \backslash \Re_{{\mathbf{C_n}}(\rm{proj} \,\Lambda)}^{{s_i}+1}$. Since, $fj_{_i}=0$ then $d_l(f) < \infty$.

$(b)$ Since $f$ is irreducible and ${\rm Coker}\,f \in \mathbf{C_n}(\rm{proj} \,\Lambda)\backslash  \{0\}$ by
 Proposition \ref{Nu-Conucleo}, $(a)$, and   Theorem \ref{Nu yconu indesc} we have that ${\rm Coker}\,f $ is indecomposable.  The proof of $(b)$
follows  with similar arguments as Statement $(a)$.
\end{proof}

\subsection{A right $\mathbf{C_n}(\rm{proj}\,\Lambda)$-approximation of the kernel}

The aim of this section is to study the right $\mathbf{C_n}(\rm{proj}\,\Lambda)$-approximations of the kernel of an irreducible morphism in $\mathbf{C_n}(\rm{proj}\,\Lambda)$ in order to get information on the left degree of such an irreducible morphism.
By \cite[proof of Theorem 4.5]{BSZ:04}, for each $X \in \mathbf{C_n}(\rm{mod}\,\Lambda)$  there is a right $\mathbf{C_n}(\rm{proj}\,\Lambda)$-approximation.

\begin{defi} \label{def cn-aprox}
A  morphism $\varphi:W \rightarrow X$ in $\mathbf{C_n}(\rm{mod}\,\Lambda)$, with $W \in \mathbf{C_n}(\rm{proj}\,\Lambda)$ is a right  $\mathbf{C_n}(\rm{proj}\,\Lambda)$-approximation of $X$ if for each morphism $h:Z \rightarrow X$ with $Z \in \mathbf{C_n}(\rm{proj}\,\Lambda)$ there is a morphism $t:Z \rightarrow W$ such that $\varphi t=h$.
\end{defi}

\begin{prop} \label{Cn-aproxim en artin} Let $\Lambda$ be an artin algebra, $f:X \rightarrow Y$ an irreducible morphism in $\mathbf{C_n}(\rm{proj}\,\Lambda)$ with $X$ or $Y$ indecomposable and ${\rm Ker}\, f \neq 0$. Assume that $\varphi:W \rightarrow {\rm Ker}\, f$ is  a right  $\mathbf{C_n}(\rm{proj}\,\Lambda)$-approximation of ${\rm Ker}\,f$  such that $j \varphi \in \Re_{\mathbf{C_n}(\rm{proj}\,\Lambda)}^m (W,  X) \backslash
\Re_{\mathbf{C_n}(\rm{proj}\,\Lambda)}^{m+1}(W,  X)$ for some $m$, where $j: {\rm Ker}\, f \rightarrow X$ is the inclusion morphism in $\mathbf{C_n}(\rm{mod}\,\Lambda)$. Then, $d_l(f) < \infty$.
\end{prop}

\begin{proof} Let $\varphi:W \rightarrow {\rm Ker}\, f$ be a right $\mathbf{C_n}(\rm{proj}\,\Lambda)$-approximation
of ${\rm Ker}\,f$ satisfying the hypothesis.
Since $f j \varphi  =0$ then $d_l(f) < \infty$.
\end{proof}

As an immediate consequence we get the following corollary.

\begin{coro} \label{cn-proy aprox en trf} Let $\mathbf{C_n}(\rm{proj}\,\Lambda)$ be of finite type with $\Lambda$ an artin algebra.
If $f:X \rightarrow Y$ an irreducible morphism in $\mathbf{C_n}(\rm{proj}\,\Lambda)$ with $X$ or $Y$ indecomposable  and if ${\rm Ker}\, f \neq 0$, then
$d_l(f) < \infty$.
\end{coro}

\begin{proof} Let  $f:X \rightarrow Y$ be irreducible in $\mathbf{C_n}(\rm{proj}\,\Lambda)$ with ${\rm Ker}\, f \neq 0$ and
$j: {\rm Ker}\, f \rightarrow X$ in $\mathbf{C_n}(\rm{mod}\,\Lambda)$ the inclusion morphism.

Assume that ${\rm Ker}\, f \in \mathbf{C_n}(\rm{proj}\,\Lambda)$. Then, $j  \in \Re_{\mathbf{C_n}(\rm{proj}\,\Lambda)}^m( {\rm Ker}\, f,X))\backslash
\Re_{\mathbf{C_n}(\rm{proj}\,\Lambda)}^{m+1}( {\rm Ker}\, f,X)$ for some positive integer $m$. Since $f j =0$ then $d_l(f) < \infty$.

On the other hand, if ${\rm Ker}\, f \notin \mathbf{C_n}(\rm{proj}\,\Lambda)$ then there exists $\varphi:W \rightarrow {\rm Ker}\,f$ a right $\mathbf{C_n}(\rm{proj}\,\Lambda)$-approximation of ${\rm Ker}\,f$. Then, $j \varphi  \in \Re_{\mathbf{C_n}(\rm{proj}\,\Lambda)}^m( W ,X))\backslash
\Re_{\mathbf{C_n}(\rm{proj}\,\Lambda)}^{m+1}( W,X)$, for some positive integer $m$. By Proposition  \ref{Cn-aproxim en artin}, $d_l(f) < \infty$.
\end{proof}

\begin{prop} \label{tipo finito}
Let $\Lambda$ be an artin algebra and $ \mathbf{C_n}(\rm{proj} \,\Lambda)$ of finite type. Let $f:X \rightarrow Y$ be an irreducible morphism in
$\mathbf{C_n}(\rm{proj} \,\Lambda)$  with $X$ or $Y$ indecomposable and  ${\rm Coker}\,f \in \mathbf{C_n}(\rm{proj} \,\Lambda)\backslash \{0 \}$. Then $d_r(f)$ is finite.
\end{prop}
\begin{proof} By Proposition \ref{Nu-Conucleo}, $(a)$, $f$ is of type $(sec)$.  Consider  $\pi: Y \rightarrow {\rm Coker}\,f$ the cokernel morphism. Then, $\pi f = 0$, with $\pi  \neq 0$.
Hence, $\pi  \in \Re_{\mathbf{C_n}(\rm{proj} \,\Lambda)}^s(Y, {\rm Coker}\,f)  \backslash \Re_{\mathbf{C_n}(\rm{proj} \,\Lambda)}^{s+1}(Y, {\rm Coker}\,f)$, for some positive integer $s$. Then, $d_r(f) < \infty$.
\end{proof}

The next result  is related to the left degree of irreducible morphisms in a  category $\mathbf{C_n}(\rm{proj} \,\Lambda)$ of finite type with $\Lambda$ a selfinjective algebra.

\begin{prop}\label{autoinyectivas}
Let $\Lambda$ be a selfinjective algebra and $\mathbf{C_n}(\rm{proj} \,\Lambda)$ be a category  of finite type. Let  $f: X \rightarrow Y$  be an  irreducible morphism in $ \mathbf{C_n}(\rm{proj} \,\Lambda)$  with $X$ or $Y$ indecomposable.  The following conditions hold.
\begin{itemize}
\item[(a)] If $f$ is of type $(sec)$ then $d_r(f)< \infty$.
\item[(b)] If $f$ is of type $(ret)$ then $d_l(f)<\infty$.
\item[(c)] If $f$ is of type $(ret-irred-sec)$ then $d_l(f)<\infty$ and $d_r(f)< \infty$.
\end{itemize}
\end{prop}

\begin{proof} From  Proposition \ref{tipo finito} and Corollary \ref{cn-proy aprox en trf}  we get Statements $(a)$ and $(b)$, respectively.

$(c).$ Let $f=\{f^j\}_{j=1}^n:X \rightarrow Y$ be of type $(ret-irred-sec)$ where $f^i:X^i \rightarrow Y^i$ is the irreducible morphism in $\rm{proj}\,\Lambda$.  Then, ${\rm Ker}\,f \notin \mathbf{C_n}(\rm{proj} \,\Lambda)$ and ${\rm Coker}\,f \notin \mathbf{C_n}(\rm{proj} \,\Lambda)$, since by  \cite[Proposition 5.2]{BS:10} we have that  ${\rm Ker}\,f^i \notin \rm{proj} \,\Lambda$ and ${\rm Coker}\,f^i \notin \rm{proj} \,\Lambda$. By Corollary \ref{cn-proy aprox en trf}, $d_l(f) <\infty$.

Let  $I$ be an indecomposable direct summand of $I_{_0}({\rm Coker}\,f^i)$, where $I_{_0}({\rm Coker}\,f^i)$ is the minimal injective envelope of ${\rm Coker}\,f^i$.
Then, there is a non-zero morphism $\pi^i:Y^i \rightarrow I$  that factors through ${\rm Coker}\,f^i$.
Consider $\pi=\{0,\cdots,\pi^i,\cdots,0\}:Y \rightarrow S_i(I)$ as follows:

$$
\def\objectstyle{\scriptstyle}
\def\labelstyle{\scriptstyle}
\xymatrix @!0 @R=10mm @C=12mm{
 X :\ar[d]^f &  X^1\ar[r] \ar[d]&   \cdots \ar[r]& X^{i-1}\ar[r]\ar[d]   & X^i\ar[r]\ar[d]^{f^i}  & X^{i+1}\ar[r]\ar[d]  & \cdots \ar[r]& X^n\ar[d]
 \\Y \ar[d]^{\pi} & Y^1\ar[r] \ar[d] &  \cdots \ar[r] &Y^{i-1}\ar[r]\ar[d] & Y^i \ar[r]\ar[d]^{\pi^i}& Y^{i+1}\ar[r]\ar[d]  & \cdots \ar[r] & Y^n \ar[d]\\
 S_i(I) &  0\ar[r] &  \cdots \ar[r] & 0\ar[r]  & I \ar[r]& 0\ar[r] &  \cdots \ar[r] & 0}
$$

\noindent Then $\pi f = 0$, because $\pi$ factors through ${\rm Coker}\,f$. Hence $\pi \in \Re_{\mathbf{C_n}}^k (Y, S_i(I)) \backslash \Re_{\mathbf{C_n}}^{k+1}(Y, S_i(I))$, for some positive integer  $k$, proving that $d_r(f) < \infty$.
\end{proof}

Now, we show an example of how to compute the degree of an irreducible morphism in $\mathbf{C_2}(\rm{proj} \,\Lambda)$ when $\Lambda$ is a selfinjective algebra.

\begin{ej}\label{ej cn-aprox}
{\em Consider $\Lambda$ the  algebra given by the quiver}

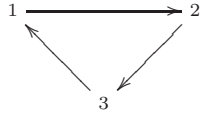
\tiny
\begin{displaymath}
\xymatrix {1 \ar[rr] & & 2\ar[ld] \\  & 3\ar[lu]&  }
\end{displaymath}
\normalsize
\vspace{.1in}

\noindent {\em with ${\rm rad}^2 \,\Lambda = 0.$ The Auslander-Reiten quiver of $\mathbf{C_2}(\rm{proj} \,\Lambda)$ is the following:}

{\tiny
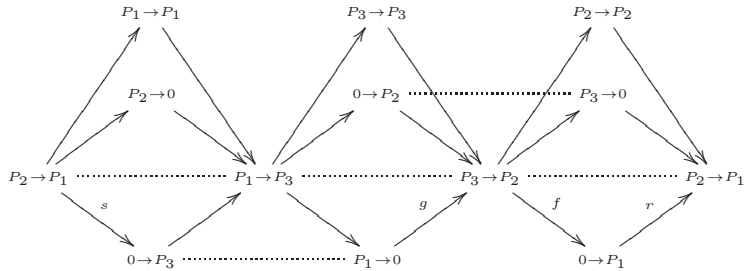
\begin{displaymath}
\def\objectstyle{\scriptstyle}
\def\labelstyle{\scriptstyle}
\xymatrix @!0 @R=1.1cm  @C=1.5cm  { &P_1 \rightarrow P_1\ar[rdd]
& & P_3 \rightarrow P_3 \ar[rdd]& &
P_2 \rightarrow  P_2\ar[rdd] &  \\
& P_2 \rightarrow 0 \ar[rd]& & 0 \rightarrow P_2\ar[rd] \ar@{.}[rr]& & P_3 \rightarrow  0 \ar[rd]& \\
P_2 \rightarrow P_1 \ar[ruu]\ar[ru]\ar[rd]^s \ar@{.}[rr] & & P_1 \rightarrow
P_3  \ar[ruu]\ar[ru]\ar[rd]\ar@{.}[rr]  &
& P_3 \rightarrow P_2 \ar[ruu]\ar[ru]\ar[rd]^f \ar@{.}[rr] & & P_2 \rightarrow P_1\\
 & 0 \rightarrow P_3\ar[ru]\ar@{.}[rr] & & P_1 \rightarrow 0 \ar[ru]^g & & 0   \rightarrow P_1\ar[ru]^r &
}
\end{displaymath}}
\vspace{.1in}

\noindent {\em The kernel of the irreducible morphism  $f: (P_3 \rightarrow P_2)\rightarrow (0\rightarrow P_1)$ is the complex $(P_3 \rightarrow S_3)$ which does not belong to $\mathbf{C_2}(\rm{proj} \,\Lambda)$. A right $\mathbf{C_2}(\rm{proj} \,\Lambda)$-approximation of such a kernel is the complex $W=(P_1 \rightarrow 0) \oplus (P_3 \rightarrow P_3)$. Then $d_l(f)=1$ since $g:(P_1 \rightarrow 0)\rightarrow (P_3 \rightarrow P_2) \in \Re_{\mathbf{C_2}(\rm{proj} \,\Lambda)} \backslash \Re^2_{\mathbf{C_2}(\rm{proj} \,\Lambda)} $ and $g f=0$.}
\end{ej}

The following result allows us to compute the left degree of certain irreducible morphisms taking into account the right $\mathbf{C_n}(\rm{proj}\,\Lambda)$-approximation  of the kernel. More precisely, we prove Theorem A.

\begin{thm} \label{Cn-aproxim} Let $\Lambda$ be an artin algebra and $\Gamma$  a generalized standard component of $\Gamma_{\mathbf{C_n}(\rm{proj}\,\Lambda)}$ with length. Let $f:X \rightarrow Y$ be an irreducible morphism with $X \in \Gamma$, $\varphi: W \rightarrow {\rm Ker}\,f$ a right
$\mathbf{C_n}(\rm{proj}\,\Lambda)$-approximation  of ${\rm Ker}\,f$ and $j:{\rm Ker}\, f \rightarrow X$ the inclusion morphism in $\mathbf{C_n}(\rm{mod}\,\Lambda)$. The following conditions are equivalent.
\begin{itemize}
\item [(a)] $d_l(f) = m$.
\item [(b)] $j \varphi \in \Re_{\mathbf{C_n}}^m(W,X)\backslash
\Re_{\mathbf{C_n}}^{m+1}(W,X)$. Moreover, there exists an indecomposable direct summand $Z$  of $W$ and a morphism $h \in   \Re_{\mathbf{C_n}(\rm{proj}\,\Lambda)}^m(Z,X)\backslash
\Re_{\mathbf{C_n}(\rm{proj}\,\Lambda)}^{m+1}(Z,X)$ such that $fh=0$.
 \end{itemize}
\end{thm}

\begin{proof} $ (a) \Rightarrow (b).$  Assume $d_l(f)= m$.  By  Theorem \ref{first}, there exist $Z \in \Gamma$,
a positive integer $m$  and a morphism $h$ such that  $h \in \Re_{\mathbf{C_n}(\rm{proj}\,\Lambda)}^m
(Z,X)\backslash  \Re_{\mathbf{C_n}(\rm{proj}\,\Lambda)}^{m+1}(Z,X)$ with $fh=0$.
Hence, $h$ factors through ${\rm Ker}\, f$ in $\mathbf{C_n}(\rm{mod}\,\Lambda)$. Hence, there is a morphism $\alpha:Z \rightarrow
{\rm Ker}\, f$  in $\mathbf{C_n}(\rm{mod}\,\Lambda)$ such that $j \alpha = h$. Then, there exists
$\alpha': Z \rightarrow W$ such that  $\varphi \alpha' =
\alpha$. We illustrate the situation with the following diagram:

$$
\xymatrix{ {\rm Ker} \,f \ar[r]^{j} & X \ar[r]^f & Y
\\W \ar[u]_{\varphi} & Z \ar@{.>}[l]^{\alpha '} \ar[u]^{h} \ar@{.>}[lu]_{\alpha} \ar@{.>}[ru]^0}$$
\vspace{.1in}

\noindent Note that $h = (j \varphi)\alpha' \in \Re_{\mathbf{C_n}}^m
(Z,X)\backslash \Re_{\mathbf{C_n}}^{m+1}(Z,X)$. Since $\Gamma$ is generalized standard  there exists $t \leq m$ such that $j \varphi \in
\Re_{\mathbf{C_n}(\rm{proj}\,\Lambda)}^t(W, X)\backslash \Re_{\mathbf{C_n}(\rm{proj}\,\Lambda)}^{t+1}(W, X)$
If $t<m$, since $f j \varphi =0$ then  $d_l(f) < m$ contradicting the hypothesis. Then, $t=m$ and $\alpha'$ is
a section, proving that $Z$ is a direct summand of $W$.

$(b) \Rightarrow (a)$. Let $\varphi: W
\rightarrow {\rm Ker}\, f$  be a right  $\mathbf{C_n}(\rm{proj}
\,\Lambda)$-approximation of ${\rm Ker}\,f$ such that $j \varphi \in \Re_{\mathbf{C_n}(\rm{proj}\,\Lambda)}^m(W,X) \backslash
\Re_{\mathbf{C_n}(\rm{proj}\,\Lambda)}^{m+1}(W,X)$. By Theorem \ref{Cn-aproxim en artin}, we know that  $d_l(f) \leq m$. Assume that $d_l(f) = s < m$.
By $(a) \Rightarrow (b)$ we infer that $j \varphi \in \Re_{\mathbf{C_n}(\rm{proj}\,\Lambda)}^s(W,X) \backslash
\Re_{\mathbf{C_n}(\rm{proj}\,\Lambda)}^{s+1}(W,X)$, getting a contradiction. Then,  $d_l(f)= m.$
\end{proof}

As an immediate  consequence of the above theorem we get the following  result.

\begin{coro} \label{car ker de here} Let $\Lambda$ be an artin algebra and $\Gamma$  a generalized standard component of $\Gamma_{\mathbf{C_n}(\rm{proj}\,\Lambda)}$ with length. Let $f:X \rightarrow Y$ be an irreducible morphism with $X \in \Gamma$ and  ${\rm Ker}\,f \in \mathbf{C_n}(\rm{proj}\,\Lambda)\backslash \{0 \}$. Then,
$d_l(f)=m$ if and only if the inclusion morphism $j \in \Re_{\mathbf{C_n}(\rm{proj}\,\Lambda)}^m ({\rm Ker}\,f , X ) \backslash
\Re_{\mathbf{C_n}(\rm{proj}\,\Lambda)}^{m+1}({\rm Ker}\,f , X )$. Moreover, if $d_l(f)=m$ there exists $h \in \Re_{\mathbf{C_n}(\rm{proj}\,\Lambda)}^m(Z,X) \backslash
\Re_{\mathbf{C_n}(\rm{proj}\,\Lambda)}^{m+1}(Z,X)$ such that $fh=0$, where $Z$ is an indecomposable direct summand of ${\rm Ker}\,f$.
\end{coro}

We end up this section by showing an example of an irreducible morphism with finite right degree  and finite left degree. We observe that such a situation  does not hold in
a category  $\rm{mod} \,\Lambda$ with $\Lambda$ a finite dimensional $k$ algebra over an algebraically closed field. This result holds as a consequence of  \cite[Theorem C]{CLT}.

\begin{ej}
{\em Consider the algebra $\Lambda$ stated in Example \ref{ej cn-aprox}. There we proved that the irreducible morphism  $f: (P_3 \rightarrow P_2)\rightarrow (0\rightarrow P_1)$ in $\mathbf{C_2}(\rm{proj} \,\Lambda)$ has  left degree one.}

{\em On the other hand, if we consider the irreducible morphisms $r : (0 \rightarrow P_1)\rightarrow (P_2 \rightarrow P_1)$ and  $s: (P_2 \rightarrow P_1)\rightarrow (0 \rightarrow P_3)$ in $\mathbf{C_2}(\rm{ proj} \,\Lambda)$ by \cite[Theorem 6.12]{CPS1} we have that $$sr
 \in \Re^{2}_{\mathbf{C_2}(\rm{proj} \,\Lambda)} ((0 \rightarrow P_1), (0 \rightarrow P_3)) \backslash \Re^3_{\mathbf{C_2}(\rm{proj} \,\Lambda)} ((0 \rightarrow P_1), (0 \rightarrow P_3)).$$ Moreover, $sr f=0$. Therefore, $d_r(f)< \infty$.}
\end{ej}

\section{Degrees of irreducible morphisms in $\mathbf{C_n}({\rm proj} \,{ H})$, with $H$ a hereditary algebra}

In \cite{CLT}, the authors gave a characterization of when the left degree of an irreducible morphism is finite for a  category of
finitely generated modules over a finite dimensional algebra over an algebraically closed field. To achieve such a result they used covering techniques. More precisely, they reduce the study of the degree of an irreducible morphism in
a component to the study of the degree of an irreducible morphism in a convenient covering.

Next, we shall use a similar technique to characterize when the left (right) degree of an irreducible morphism is finite in a category $\mathbf{C_n}({\rm proj} \,H)$, for $H$ a finite dimensional hereditary algebra.

\begin{sinob} {\bf Some preliminaries on covering techniques.} Let $\Gamma$ be a translation quiver, that is, $\Gamma$ is a quiver with no loops (but with possibly parallel arrows), endowed with two distinguished subsets of vertices, the elements of which are called non projective and non injective objects, respectively; and endowed with a bijection $\tau:x \rightarrow \tau x$ (the translation) from the set of non-projective to the set of non-injective; such that for every vertices $x$, $y$ with $x$ non-projective, there is a bijection $\alpha \rightarrow \sigma\alpha$ from the set of arrows $y \rightarrow x$ to the set of arrows $\tau x \rightarrow y$. All translation quivers are assumed to be locally finite. The subquiver of $\Gamma$ formed by the arrows starting in $\tau x$ and the arrows arriving in $x$ is called the {\bf mesh} ending at $x$. We write $k(\Gamma)$ for the {\bf mesh-category} of $\Gamma$, which is the factor category of the path category $k \Gamma$ by the ideal generated by the morphisms $\sum_{\alpha:.\rightarrow x} \alpha \sigma\alpha$ where $\alpha$ goes through the arrows arriving at $x$, for a given non-projective vertex $x$.

We denote by ${\rm rad}\, k(\Gamma)$ the ideal of $k(\Gamma)$ generated by
 $\{{\overline{\alpha}} \mid \alpha \in \Gamma_1 \}$ where $\Gamma_1$ is the set of arrows in $k\Gamma$;  by
${\rm rad}^0 k(\Gamma)= k(\Gamma) $, ${\rm rad}^1 k(\Gamma)={\rm rad}\, k(\Gamma)$ and
${\rm rad}^{m+1} k(\Gamma)= {\rm rad} k(\Gamma){\rm rad}^m k(\Gamma)(= {\rm rad}^m k(\Gamma){\rm rad}
k(\Gamma)$). The radical ${\rm rad}\, k(\Gamma)$ satisfies the following useful result proven in \cite[Proposition 2.1]{CT}, which is a more general version of Proposition \ref{R}.

\begin{lem} Let $\Gamma$ be a translation quiver with length.  Assume there is a path from $x$ to $y$ in $\Gamma$ of length $m$. Then, the following conditions hold.
\begin{enumerate}
\item[(a)] $k(\Gamma)(x,y)= {\rm rad}(x,y)= \cdots ={\rm rad}^m (x,y)$,
\item[(b)] ${\rm rad}^{i}(x,y)=0$ if $i> m$.
\end{enumerate}
\end{lem}
\end{sinob}

\begin{sinob}  A covering of translation quivers is a morphism
$p:\Gamma' \rightarrow \Gamma$ of quivers such that:\begin{enumerate}
\item[(a)] $\Gamma'$ is a translation quiver.
\item[(b)] A vertex $x \in \Gamma'$  is projective (or injective, respectively) if and only if so is $p(x)$.
\item[(c)] $p$ commutes with the translations in $\Gamma$ and $\Gamma'$, that is, $p(\tau'x)= \tau p(x)$.
\item[(d)] For every vertex $x \in \Gamma'$, the map $\alpha \rightarrow p(\alpha)$ induces a bijection from the set of arrows in $\Gamma'$ starting from $x$ (or ending at $x$) the set of arrows in $\Gamma$ starting from $p(x)$ (or ending at $p(x)$, respectively).
\end{enumerate}

In \cite{CLT}, the authors defined an equivalence relation $\sim$  on the set of unoriented paths in $\Gamma$ satisfying the following properties:

\begin{enumerate}
\item[(i)] If  $\alpha:x \rightarrow y$ is an arrow in  $\Gamma$ then
$\alpha^{-1}\alpha \sim e_x$, $\alpha\alpha^{-1} \sim
e_y$ (where $e_x$ denotes the stationary path at $x$ of length $0$).

\item[(ii)] If $x$ is a non-projective vertex and the mesh in $\Gamma$ ending at $x$ has the form

{\tiny
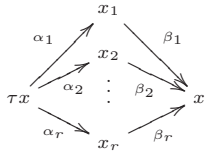
$$
\xymatrix @!0 @C=12mm @R=6mm { & x_1 \ar[rdd]^{\beta_1} & \\& x_2 \ar[rd]_{\beta_2}
\\\tau x \ar[ruu]^{\alpha_1}\ar[ru]_{\alpha_2} \ar[rd]_{\alpha_r} & \vdots
&x
\\  & x_r \ar[ru]_{\beta_r}  }
$$ }

\noindent then $\beta_i \alpha_i \sim \beta_j \alpha_j $, for every $i,j =1,\cdots,r$.

\item[(iii)] If $v \sim v'$ then $wv \sim wv'$ y $vu \sim
v' u$, where the compositions are defined.

\item[(iv)] If $\alpha$, $\beta$ are arrows in $\Gamma$ with the same source and the same  target then
 $\alpha \sim \beta$.
\end{enumerate}

Applying the construction given in \cite[(1.3)]{BG:82} and taking into account the above equivalence relation the authors obtained a covering which they called the generic covering and they denoted it by $\pi: \widetilde{\Gamma} \rightarrow \Gamma$.
By \cite[Proposition 1.2]{CLT}, we know that the generic covering of $\Gamma$ is a quiver with
length. For more details on the generic covering we refer the reader to \cite{CLT}.

We recall that a $k$-lineal functor  $F: k({\widetilde{\Gamma}}) \rightarrow {\rm{add} (\rm{ind}\,\Gamma)}$ is called  well-behaved if it  satisfies the following conditions for every vertex  $x \in
{\widetilde{\Gamma}}$: (a) $Fx= \pi x$;
(b) If $\alpha_1: x\rightarrow x_1, \cdots,\alpha_r:x
\rightarrow x_r$ are all arrows in $\widetilde{\Gamma}$ starting in $x$ then
$[F(\overline{\alpha}_1),\cdots,F(\overline{\alpha}_r)]^t:Fx
\rightarrow F x_1\oplus\cdots \oplus F x_r$ is  a minimal left almost split morphism in $\mathbf{C_n}({\rm proj}\,H)$;
(c) If $\alpha_1: x_1\rightarrow x, \cdots,\alpha_r:x_r
\rightarrow x$ are all arrows in $\widetilde{\Gamma}$ ending in
 $x$ then
$[F(\overline{\alpha}_1),\cdots,F(\overline{\alpha}_r)]:
 F x_1\oplus\cdots \oplus F x_r  \rightarrow Fx$ is a minimal right almost split morphism in $\mathbf{C_n}({\rm proj}\,H)$.
Note that these conditions imply that $F$ maps meshes of $\widetilde{\Gamma}$ to almost split sequences in
$\mathbf{C_n}({\rm proj}\,H)$.

Observe that \cite[Proposition 2.6]{CLT} and \cite[Lemma 2.7]{CLT} can be adapted to $\mathbf{C_n}({\rm proj}\,H)$ with similar proofs.
For the convenience of the reader we state both results and we refer the reader to \cite{CLT} for their proofs.

\begin{prop}\emph{(\cite[Proposition 2.6]{CLT}) \label{F(alpha)=f}} Let $f=[f_1,\cdots,f_r]^t:X \rightarrow \oplus_{i=1}^r X_i$ be an irreducible morphism
in $\mathbf{C_n}({\rm proj}\,{H})$ with $X, X_1, \cdots, X_r \in
\Gamma$. Let $x \in
\pi^{-1}(X)$ and $\alpha_i :x \rightarrow x_i$ be an arrow in
$\widetilde{\Gamma}$ such that $\pi x_i =X_i$, $\pi(\alpha_i): X \rightarrow X_i$, for $i=1,\cdots,r$. Then, there exists a well-behaved
functor  $F:k(\widetilde{\Gamma})
\rightarrow {\rm ind}\;\Gamma$ such that
$F(\overline{\alpha}_i)=f_i$, for every $i$.
\end{prop}

\begin{lem}\emph{(\cite[Lemma 2.7]{CLT}) \label{existehi}} Let $F:k(\widetilde{\Gamma})
\rightarrow {\rm ind}\;\Gamma$ be a well behaved functor, $x$, $y$
vertices in $\widetilde{\Gamma}_0$, $m \geq 0$. Then, the following conditions hold.
\begin{enumerate}
\item[(a)] $F$ maps a morphism in ${\rm rad}^m k(\widetilde{\Gamma})(x,y)$ onto a morphism in
$\Re_{\mathbf{C_n}({\rm proj}\,H)}^m(Fx,Fy)$.
\item[(b)] Let $f \in \Re_{\mathbf{C_n}({\rm proj}\,H)}^{m+1}(Fx,Fy)$ and  $\alpha_i:x \rightarrow
x_i$ , $i=1,\cdots,r$ be the arrows in $\widetilde{\Gamma}$ starting in $x$. Then, there exists  $h_i \in
\Re_{\mathbf{C_n}({\rm proj}\,H)}^m(Fx_i,Fy)$, for every $i$, such that, $f= \Sigma_{i=1}^r
h_iF(\overline{\alpha}_i)$.
\end{enumerate}
\end{lem}
\end{sinob}

\begin{rem}\label{Aut(X)en cnproyH}{\em Let $X$ be an indecomposable complex in $\mathbf{C_n}({\rm proj}\,H)$.  By
$(\ref{1.1})$ if $X$ is not injective then ${\rm Aut}_{\mathbf{C_2}({\rm proj}\,H)}(X) \simeq {\rm Aut}_{H}(M)\simeq k$ for $M$ an indecomposable $H$-module. Hence,  ${\rm Aut}_{\mathbf{C_n}({\rm proj}\,H)}(X) \simeq  k$.}

{\em  In case  $X$ is an indecomposable injective then $X=J_i(P)$ or $X=S(P)$ with $P$ an indecomposable projective module. Hence, ${\rm Aut}_{\mathbf{C_n}({\rm proj}\,H)}(X) \simeq {\rm Aut}_{H}(P) \simeq k$.}
\end{rem}

Our next result is similar to \cite[Teorema B]{CLT}, but  the proof of $(a)$ is not the same  because in $\mathbf{C_n}({\rm proj} \,H)$ an irreducible morphism  with indecomposable injective domain is not necessarily an epimorphism. We only transcribe Statement $(a)$ and  give a proof for it.

\begin{thm}\label{radicalescub} Let $F:k(\widetilde{\Gamma})
\rightarrow {\rm ind}\,\Gamma$ be a well-behaved  functor, $x$, $y$
vertices in $\widetilde{\Gamma}_0$, $m \geq 0$.
The following maps induced by $F$ are
bijective.
\begin{enumerate}
\item [(a)]$F_m: \bigoplus_{z/Fz=Fy}{\rm rad}^m k(\widetilde{\Gamma})(x,z)/ {{\rm rad}}^{m+1}
k(\widetilde{\Gamma})(x,z)\rightarrow
\Re_{\mathbf{C_n}({\rm proj}\,H)}^m(Fx,Fy)/\Re_{\mathbf{C_n}({\rm proj}\,H)}^{m+1}(Fx,Fy).$

\item [(b)]$G_n:\bigoplus_{z/Fz=Fx}{{\rm rad}}^m k(\widetilde{\Gamma})(z,y)/ {{\rm rad}}^{m+1}
k(\widetilde{\Gamma})(z,y)\rightarrow
\Re_{\mathbf{C_n}({\rm proj}\,H)}^m(Fx,Fy)/\Re_{\mathbf{C_n}({\rm proj}\,H)}^{m+1}(Fx,Fy).$
\end{enumerate}
\end{thm}

\begin{proof} $(a)$. We observe that, as in \cite[Teorema B]{CLT}, the surjective property follows by induction on $m$. We only prove that for all $m \geq 0$,  $F_m$ is injective which is equivalent to prove that for all $m \geq 0$, if
$(\phi_z)_z\in
\bigoplus_{Fz=Fy}k(\widetilde{\Gamma})(x,z)$ is such that $\Sigma_z
F(\phi_z) \in \Re_{\mathbf{C_n}({\rm proj}\,H)}^{m+1}(Fx,Fy)$ then $\phi_z \in {\rm rad}^{m+1} k(\widetilde{\Gamma})(x,z),$ for every $z$.

We  prove the result by
induction on $m \geq 0$. Consider $m=0$ and $\Sigma_z F(\phi_z) \in \Re_{\mathbf{C_n}({\rm proj}\,H)}(Fx,Fy)$. In case $Fx
\neq Fy$ then $x \neq z$ for every $z$ such that $Fz =Fy$. Therefore,
 $\phi_z \in {{\rm rad}} k(\widetilde{\Gamma})(x,z)$. In case $Fx
= Fy$  then $\phi_z \in {{\rm rad}} k(\widetilde{\Gamma})(x,z)$
if $x \neq z$ and there exists $\lambda \in k$ such that $\phi_x = \lambda
1_x$. Then, $ \lambda 1_{Fx} \in \Re_{\mathbf{C_n}({\rm proj}\,H)}(Fx,Fy)$,
and therefore, $\lambda =0$. Hence, $\phi_z \in {{\rm rad}}
k(\widetilde{\Gamma})(x,z)$, for every $z$, proving that the result is
true for $m=0$.

Now, by inductive hypothesis we have that $(\phi_z)_z\in
\bigoplus_{Fz=Fy}k(\widetilde{\Gamma})(x,z)$ is such that $\Sigma_z
F(\phi_z) \in \Re_{\mathbf{C_n}({\rm proj}\,H)}^{m+1}(Fx,Fy)$ then $\phi_z \in {\rm rad}^{m+1} k(\widetilde{\Gamma})(x,z),$ for every $z$.
Let
$(\phi_z)_z \in \bigoplus_{Fz =Fy}k(\widetilde{\Gamma})(x,z)$ be such that $\Sigma_z F(\phi_z) \in \Re_{\mathbf{C_n}({\rm proj}\,H)}^{m+2}(Fx,Fy)\subset
\Re_{\mathbf{C_n}({\rm proj}\,H)}^{m+1}(Fx,Fy).$
Then, $\phi_z \in {{\rm rad}}^{m+1} k(\widetilde{\Gamma})(x,z)$, for every $z$, since the result holds for $m$. Because $F_{m+2}$ is surjective and $\Sigma_z
F(\phi_z) \in \Re_{\mathbf{C_n}({\rm proj}\,H)}^{m+2}(Fx,Fy)$, for every $z$
 such that $Fz=Fy$ then there exists $\psi_z \in {{\rm rad}}^{m+2}
k(\widetilde{\Gamma})(x,z)$ such that $\Sigma_z F(\phi_z)=
\Sigma_z F(\psi_z)$
(modulo  $\Re_{\mathbf{C_n}({\rm proj}\,H)}^{m+3}$).
Hence,  $\Sigma_z F(\phi_z-\psi_z) \in
\Re_{\mathbf{C_n}({\rm proj}\,H)}^{m+3}(Fx,Fy)$. By Lemma \ref{existehi}, we have that
$\Sigma_z F(\phi_z-\psi_z) = \Sigma_{i=1}^r h_i
F(\overline{\alpha}_i)$ with $h_i \in
\Re_{\mathbf{C_n}({\rm proj}\,H)}^{m+2}(Fx_i,Fy)$ and where $x_i$ are the end points of all the arrows $x \rightarrow x_i$ in $\widetilde{\Gamma}$ starting in $x$, for $i=1, \dots, r$.
Moreover, since $\phi_z$ and
$\psi_z$ are in ${{\rm rad}}^{m+1} k(\widetilde{\Gamma})(x,z)
\subset {{\rm rad}} k(\widetilde{\Gamma})(x,z)$ then
$\phi_z-\psi_z = \Sigma_{i=1}^r \theta_{z,i}\overline{\alpha}_i$
with $\theta_{z,i} \in  k(\widetilde{\Gamma})(x_i,z)$ for
$i=1,\cdots r$. Hence, $\Sigma_z F(\phi_z-\psi_z)=\Sigma_z
F(\Sigma_{i=1}^r  \theta_{z,i}\overline{\alpha}_i) =
\Sigma_{i=1}^r (\Sigma_z F (\theta_{z,i})
F(\overline{\alpha}_i))$ and we  get that $\Sigma_{i=1}^r h_i
F(\overline{\alpha}_i) =\Sigma_{i=1}^r (\Sigma_z F (\theta_{z,i})
F(\overline{\alpha}_i))$. Then

\begin{equation}\label{1} \Sigma_{i=1}^r(\Sigma_z F
(\theta_{z,i})-h_i) F(\overline{\alpha}_i)=0.
\end{equation}

Now, we have two cases to analyze. First, if $x$ is an injective vertex then $Fx$ is injective in $\mathbf{C_n}({\rm proj}\,H)$.  Hence, either $Fx = J_i(P)$ for $1 \leq i \leq n-1$ or $Fx= S(P)$, with   $P$ an indecomposable projective in ${\rm mod}\,H$.

If $Fx = J_i(P)$ then $r=1$ and $F(\overline{\alpha}_1):J_i(P) \rightarrow Fx_1$ and $Fx_1=L_i(P)$, see  $(1.5)$. Since $H$ is hereditary then $L^{j}=0$ for $j \geq 2$. We have that  $(\ref{1})$ is the composition:

{\tiny
$$ \xymatrix @!0 @C=16mm {J_i(P)\ar[d]^{F(\overline{\alpha}_1)} &\cdots\ar[r] & 0 \ar[r]\ar[d] & P \ar[r]^{id}\ar[d]^{id} & P \ar[r] \ar[d]^{\lambda}& 0
\ar[r]\ar[d] & \cdots\\
               Fx_1\ar[d]^f  & \cdots\ar[r] & 0  \ar[r]\ar[d]^{0} & P \ar[r]^{\lambda}\ar[d]^{f^i}&  L^1 \ar[r] \ar[d]^{f^{i+1}} &  0
           \ar[d]^{0}  \ar[r] & \cdots
              \\ Fz & \cdots\ar[r] & Z^{i-1} \ar[r] & Z^{i} \ar[r]&  Z^{i+1} \ar[r] &  Z^{i+2}  \ar[r] & \cdots}$$
}

\noindent where $f= \Sigma_z F(\theta_{z,1})-h_1$.  Since $(\Sigma_z F
(\theta_{z,1})-h_1) F(\overline{\alpha}_1)=0$ then  $ f^{i+1}\lambda =0$. Because
$\lambda$ is irreducible in ${\rm proj}\,H$, by Corollary \ref{Cokernoproy} $(b)$,  we obtain that $f^{i+1}=0$. Then, $f^i= f^{i} id =0$.  Therefore,  $f= \Sigma_z F
(\theta_{z,1})-h_1=0$, that is, $\Sigma_z F
(\theta_{z,1})=h_1 \in \Re_{\mathbf{C_n}({\rm proj}\,H)}^{m+2}(Fx_1, Fy) \subset \Re_{\mathbf{C_n}({\rm proj}\,H)}^{m+1}(Fx_1, Fy)$. By inductive hypothesis,   $\theta_{z,1} \in {\rm rad}^{m+1}k(\widetilde{\Gamma})(x_1,z)$. Then, $\theta_{z,1}\overline{\alpha}_1 \in {\rm rad}^{m+2}k(\widetilde{\Gamma})(x,z)$. Because $\psi_z \in {\rm rad}^{m+2}k(\widetilde{\Gamma})(x,z)$ we get that $\phi_z= \psi_z + \theta_{z,1}\overline{\alpha}_1 \in {\rm rad}^{m+2}k(\widetilde{\Gamma})(x_1,z)$.

Now, if $Fx = S(P)$ then $r\geq 1$ and $(F(\overline{\alpha}_1),\cdots, F((\overline{\alpha}_r))^t:S(P) \rightarrow \oplus_{i=1}^{r} Fx_i$

{\tiny
$$
\xymatrix @!0 @C=16mm {S(P) \ar[d]^{[F(\overline{\alpha}_1),\cdots, F((\overline{\alpha}_r)]^t} & & P \ar[r]\ar[d]^{\lambda} & 0 \ar[r]\ar[d] & 0 \ar[r]
\ar[d] & \cdots
           \\  \bigoplus Fx_i \ar[d]^f  &  &  L^1 \ar[r]\ar[d]^{f^1} &   0 \ar[r] \ar[d] &  0 \ar[d]  \ar[r] & \cdots
            \\   Fz &  & Z^{1} \ar[r] & Z^{2} \ar[r]&  Z^{3} \ar[r] & \cdots      }
$$
}

\noindent where $f= [\Sigma_z F(\theta_{z,1})-h_1, \cdots,\Sigma_z F(\theta_{z,r})-h_r]$.
By Corollary \ref{Cokernoproy} $(b)$,
we have that $f^1=0$. Therefore, $f=0$. Then $\Sigma_z F
(\theta_{z,i})-h_i = 0$ for $i=1, \cdots, r$. Similarly to the previous  case we can  prove that $\phi_z \in {\rm rad}^{m+2}k(\widetilde{\Gamma})(x,z)$.

Now, assume that $x$ is not injective. Then, in  $\widetilde{\Gamma}$ we have the following situation:

{\tiny
$$
\xymatrix @!0 @C=10mm @R=10mm {\ & x_1 \ar[rd]^{\beta_1} & \\\,\,\, \,\,\,x\,\,\, \,\,\,\ar[ru]^{\alpha_1} \ar[rd]_{\alpha_r} & \vdots
& \tau^{-1} x
\\  & x_r \ar[ru]_{\beta_r}  }
$$
}

\noindent Since $F$ is a well-behaved, there is an almost split sequence in $\mathbf{C_n}({\rm proj}\,H)$ as follows

\begin{displaymath}
\xymatrix  {  0  \ar[r] & Fx \ar[rrr]^{[F((\overline{\alpha}_1)),\cdots, F((\overline{\alpha}_r))]^t } &  & & \bigoplus_{i}Fx_i  \ar[rr]^{ [F(\overline{\beta}_1),\cdots, F(\overline{\beta}_r)]} & &A_n^{-1}Fx \ar[r] &   0.}
\end{displaymath}

From $(\ref{1})$, we deduce that there exists $h \in {\rm Hom}_{\mathbf{C_n}(\rm {proj}\,H)}(A_n^{-1}Fx,Fy)$ such that $\Sigma_z F
(\theta_{z,i})-h_i =h F(\overline{\beta}_i)$, for every $i$, since $[\Sigma_z F
(\theta_{z,1})-h_1, \cdots,\Sigma_z F
(\theta_{z,r})-h_r]$ factors through the cokernel of $[F((\overline{\alpha}_1)),\cdots, F((\overline{\alpha}_r))]^t$.
Since $F_0$, $F_1$, $\cdots, F_{m-1}$, are surjective, there exists $(\chi_z)_z \in \oplus_{_{Fz=Fy}}
k(\widetilde{\Gamma} )(\tau^{-1}x,z)$ such that $h=\Sigma_zF(\chi_z)$ (modulo $\Re_{\mathbf{C_n}({\rm proj}\,H)}^m$).
Then for every $i$, the following equality holds:
$\Sigma_z F(\theta_{z,i}) = \Sigma_z F(\chi_z \overline{\beta}_i) + h_i \,\,\, ( \mbox{modulo} \,\,
  \Re_{\mathbf{C_n}({\rm proj}\,H)}^{m+1}).$
 \noindent Therefore,  $\Sigma_z F(\theta_{z,i} - \chi_z \overline{\beta}_i) \in \Re_{\mathbf{C_n}({\rm proj}\,H)}^{m+1}(Fx_i,Fz)$, because $h_i  \in \Re_{\mathbf{C_n}{\rm proj}\,H)}^{m+2}(Fx_i,Fy)$. By inductive hypothesis,  $\theta_{z,i} - \chi_z \overline{\beta}_i \in {\rm rad}^{m+1}k(\widetilde{\Gamma})(x_i,z)$, for every $i$, $z$. Then, for every $z$  we get:

 $$ \phi_z= \psi_z + \sum_i (\theta_{z,i} - \chi_z \overline{\beta}_i) \overline{\alpha}_i \in{\rm rad}^{m+2} k(\widetilde{\Gamma})(x,z).$$

\noindent  This proves that $F_{m+1}$ is injective for every integer $m$ such that $F_{m}$ is injective. Therefore, for every $m \geq 0$ the map $F_m$ is bijective.

$(b)$. The proof of  $(b)$
is  analogous to \cite[Teorema B, (a)]{CLT}.
\end{proof}

Now, we are in position to state the main result of
this section, Theorem B,  which is similar to \cite[Teorema C]{CLT} and where the proof follows with a similar argument. We refer the reader to \cite[Teorema C]{CLT} for a  proof. Observe that Theorem \ref{radicalescub} is fundamental to obtain the proof, and that we adapted such result for a category $\mathbf{C_n}({\rm proj} \,H)$ with $H$ a hereditary algebra.

\begin{thm} \label{caract grado} Let $H$ be a  finite dimensional hereditary algebra over an  algebraically closed field and $f: X \rightarrow Y$ be an irreducible morphism in $\mathbf{C_n}({\rm proj} \,H)$ with
$X$ indecomposable. Let  $\Gamma$ be an Auslander-Reiten component  of
  $\Gamma_{\mathbf{C_n}({\rm proj} \,H)}$ containing
$X$ and $m \in \mathbb{N}$.
\begin{enumerate}
  \item [(a)] If $d_l(f) =m$, then there exists $Z \in \Gamma$ and $h
  \in \Re_{\mathbf{C_n}({\rm proj} \,H)}^{m} (Z,X) \backslash \Re_{\mathbf{C_n}({\rm proj} \,H)}^{m
  +1}(Z,X)$ such that $fh =0$.

  \item [(b)] If $d_r(f) =m$, then there exists $Z \in \Gamma$ and $h
  \in \Re_{\mathbf{C_n}({\rm proj} \,H)}^{m} (Y,Z) \backslash \Re_{\mathbf{C_n}({\rm proj} \,H)}^{m
  +1}(Y,Z)$ such that $hf =0$.
\end{enumerate}
\end{thm}

Following Proposition \ref{Ker en Cn(H)} and  Theorem  \ref{Nu yconu indesc} if the kernel of an irreducible morphism $f$ is non-zero then we know that it lies in
$\mathbf{C_n}({\rm proj} \,H)$ and furthermore it is indecomposable.

\begin{coro}\label{ker en gamma} Let $H$ be a finite dimensional hereditary algebra over an algebraically closed field  and $\Gamma \subset \Gamma_{\mathbf{C_n}({\rm proj} \,H)}$. If $f:X \rightarrow Y$ is an irreducible morphism in ${\mathbf{C_n}({\rm proj} \,H)}$ with $X \in \Gamma$ such that $d_l(f)=m$ then the inclusion morphism $j:{\rm Ker}\,f \rightarrow X$ is such that
$$j \in \Re_{\mathbf{C_n}({\rm proj} \,H)}^m ({\rm Ker}\,f , X)\backslash \Re_{\mathbf{C_n}({\rm proj} \,H)}^{m+1}({\rm Ker}\,f , X).$$
\end{coro}

\begin{rem}\label{conucleo} We observe that Corollary \ref{ker en gamma} has a dual version when $Y \in \Gamma$  and  $d_r(f)=m$ if ${\rm Coker}\,f \in \mathbf{C_n}({\rm proj} \,H)$. In fact, in that case the cokernel morphism $\Pi$ is such that  $\Pi:Y \rightarrow {\rm Coker}\,f \in \Re_{\mathbf{C_n}({\rm proj} \,H)}^m (Y,{\rm Coker}\,f)\backslash \Re_{\mathbf{C_n}({\rm proj} \,H)}^{m+1}(Y,{\rm Coker}\,f).$
\end{rem}

\subsection{The type of $\mathbf{C_n}({\rm proj} \,H)$}
In \cite[Theorem A]{CLT} the authors gave a characterization of when a finite dimensional $k$-algebra over an algebraically closed field
is representation-finite taking into account the degree of some particular irreducible morphisms. In this section we are interested to give an equivalent  result for $\mathbf{C_n}({\rm proj} \,H)$.

We recall that $\mathbf{C_n}({\rm proj} \,H)$ is  of finite type if and only if $H$ is a path algebra
of a Dynkin quiver.

  We prove Theorem C.

\begin{thm}\label{tipofinito}
Let $H$ be a finite dimensional $k$-algebra over an algebraically closed field. The following conditions are equivalents.
\begin{enumerate}
\item [(a)] The category $\mathbf{C_n}({\rm proj} \,H)$ is of finite type.
\item [ (b)]  There exists a simple projective $H$-module $P$ such that the irreducible morphism ${\rho}_{n-1}^n:T(P) {\rightarrow} J_{n-1}(P)$
         in $\mathbf{C_n}({\rm proj} \,H)$  has finite right degree.
\item [(c)] For all simple projective $H$-module $Q$,  the irreducible morphism
         ${\rho}_{n-1}^n:  T(Q) {\rightarrow} J_{n-1}(Q)$ in $\mathbf{C_n}({\rm proj} \,H)$ has finite right degree.
\item [(d)] There exists a  projective $H$-module $P$ with $\nu(P)$ a  simple injective $H$-module such that the irreducible morphism ${\lambda}^n_1: J_1(P){\rightarrow} S(P)$  in $\mathbf{C_n}({\rm proj} \,H)$ has finite left degree.
\item [(e)]  For all projective $H$-module $Q$ with $\nu(Q)$ a simple injective $H$-module, the irreducible morphism
        ${\lambda}^n_1: J_1(Q) {\rightarrow} S(Q)$  in
        $\mathbf{C_n}({\rm proj} \,H)$ has finite left degree.
\item [(f)] All irreducible morphisms $f:X \rightarrow Y$ with $X$ or $Y$ indecomposable complexes in $\mathbf{C_n}{(\rm proj} \,H)$ of type $(sec)$
        are such that $d_r(f) < \infty$.
\item [(g)] All irreducible morphisms $f:X \rightarrow Y$ with $X$ or $Y$ indecomposable complexes  in  $\mathbf{C_n}({\rm proj} \,H)$ of type $(ret)$ are such that $d_l(f)< \infty$.
      \end{enumerate}
\end{thm}
\begin{proof} The implications $(g) \Rightarrow (e) \Rightarrow (d)$ and $(f) \Rightarrow (c) \Rightarrow (b)$ are immediate.

$(b) \Rightarrow (a)$. Let
 $\rho^n_{n-1}: T(P) {\rightarrow} J_{n-1}(P)$ be an irreducible morphism in $\mathbf{C_n}({\rm proj} \,H)$ with $d_r({\rho}^n_{n-1})=m <\infty$ and where $P$ is a simple projective  $H$-module. Let $\Gamma$ be a component of $\Gamma_{\mathbf{C_n}({\rm proj} \,H)}$ that contains the projective $J_{n-1}(P)$. By  Remark \ref{conucleo}  we have that ${\rm Coker}\,\rho^n_{n-1}= S_{n-1}(P) \in \Gamma$.  Consider the embedding $\varphi_{n-2}:\mathbf{C_2}({\rm proj} \,H)\rightarrow \mathbf{C_n}({\rm proj} \,H) $ defined by
$\varphi_{n-2}(X)= X[n-2]$, see (\ref{1.1}). Consider the morphism $\rho^2_1:(0 \rightarrow P)\rightarrow (P \rightarrow P)$ in $\mathbf{C_2}({\rm proj} \,H)$, and ${\rm Coker}\,\rho^2_1=(P \rightarrow 0)$. Then, $\varphi_{n-2}(\rho^2_1)=\rho^n_{n-1}$ and $\varphi_{n-2}({\rm Coker}\,\rho^2_1)= {\rm Coker}\,\rho^n_{n-1}$. By the shape of the Auslander-Reiten quiver of $\mathbf{C_n}({\rm proj} \,H)$, which are copies of the Auslander-Reiten quiver of $\mathbf{C_2}({\rm proj} \,H)$, the complexes
$(0 \rightarrow P)$, $(P \rightarrow P)$ and $(P \rightarrow 0)$ are in the same component $\Gamma'$ of $\Gamma_{\mathbf{C_2}({\rm proj} \,H)}$.

To prove that $\mathbf{C_n}({\rm proj} \,H)$  of finite type is equivalent to showing  that $\mathbf{C_2}({\rm proj} \,H)$ is also of finite type.

Consider $j$ a vertex in $Q_H$. Then, there is a walk of irreducible morphisms between indecomposable projective $H$-modules from $P$ to $P_j$ in $\Gamma_H$. By \cite[Proposition 3.2]{CPS1} there is a walk of irreducible morphisms from  $(0 \rightarrow P)$ to $(0 \rightarrow P_j)$ in $\Gamma'$. Hence, for each vertex $j \in Q_H$ we get that all indecomposable projective objects  of the form $(0 \rightarrow P_j)$ are in the component $\Gamma'$.

Since   ${\rm Coker}\,\rho^2_1=(P\rightarrow 0) \in \Gamma'$, we consider for each vertex $j$ in $Q_H$ a walk of irreducible morphisms between indecomposable projective $H$-modules from $P$ to $P_j$ in $\Gamma_H$. By  \cite[Proposition 3.2]{CPS1},  there is a walk  of irreducible morphisms in $\Gamma'$  from $(P \rightarrow 0)$ to $(P_j  \rightarrow 0)$. Then, for each vertex $j \in Q_H$ the injective complexes $(P_j \rightarrow 0)$ belong to the same component that the complex $(0 \rightarrow P)$.
By the equivalence between $\mathbf{\overline{C}_2}({\rm proj} \,H)$ and ${\rm mod}\,H$ given in \cite[Proposition 3.3]{B:04},
a projective of the form $(0 \rightarrow P_j)$ is in correspondence with $P_j$  and the complex $A_2(P_j \rightarrow 0)$ with the injective $I_j$. Then, the projective and the injective $H$-modules are all in the same component of $\Gamma_H$. By  \cite[Proposition 1.3]{ARS:95}
we have  that ${\rm mod}\,H$ is representation-finite and hence $\mathbf{C_2}({\rm proj} \,H)$ is of finite type.
Similarly, we get that $(d) \Rightarrow (a).$

$(a)\Rightarrow (f).$ Since $f$ is of type $(sec)$ then by Proposition \ref{Nu-Conucleo}, $(a)$ we have that ${\rm Coker}\,f \in \mathbf{C_n}({\rm proj}\,H)$. Consider
$\pi:Y \rightarrow  {\rm Coker}\,f$, the cokernel morphism of $f$. Since  $\mathbf{C_n}({\rm proj} \,H)$ is of finite type then, $\pi \in \Re^m_{\mathbf{C_n}({\rm proj} \,H)} (Y,{\rm Coker}\,f)\backslash \Re^{m + 1}_{\mathbf{C_n}({\rm proj} \,H)}(Y,{\rm Coker}\,f)$ for some $m$. Hence $d_r(f) < \infty$ because $\pi f=0$.
With a similar argument we  prove $(a)\Rightarrow (g).$
\end{proof}

\section{Nilpotency index of the radical of $\mathbf{C_n}({\rm proj} \,H)$}

In this section,  we are interested in studying  which is the less positive integer $r$ such that $\Re^{r}({\mathbf{C_n}({\rm proj} \,H)})=0$.
To achieve this, we first study the radical of the category $\mathbf{C_2}({\rm proj} \,H)$. Then, we infer the results for any $n \geq 4$ by knowing the corresponding results for   $\mathbf{C_3}({\rm proj} \,H)$ which are slightly different from $\mathbf{C_2}({\rm proj} \,H)$.

For $i$  a sink  and  $j$  a source in $Q_H$, respectively, and  $Q_H$ a Dynkin diagram,
we denote by $\rho_k^{n,i}:R_k(P_i) \rightarrow J_k(P_i)$ and by $\lambda_k^{n,j}:J_k(P_j)  \rightarrow L_k(P_j)$ the  irreducible morphisms in $\mathbf{C_n}({\rm proj} \,H)$ and $k=1,\cdots, n-1$. Note that these morphisms are the same as the ones defined in $(1.5)$ in $\mathbf{C_n}({\rm proj} \,\Lambda)$ for any $n \geq 2$.

\begin{rem}\label{r_H} \emph{ By \cite[Theorem 2.6]{C} the nilpotency index of the radical of $\rm {mod}\,\Lambda$, for $\Lambda$ a finite dimensional algebra over an algebraically closed field of finite representation type, is equal to one plus the length of any path from $P_i$ to
 $I_i$ going through the simple $S_i$.}
 \end{rem}


First, note that for any hereditary algebra of  Dynkin type,
we have that $\Gamma_{H}$ is a component with length. Indeed the result follows from the characterization stated in Theorem \ref{Gamasc}, which was given in \cite{BG:82} in the proof of Proposition 6.

We recall that a translation quiver is say to be simply connected if it is connected  and the fundamental group \mbox{$\pi _1(\Gamma,x)=1$} for some vertex $x\in \Gamma$.

\begin{thm}\label{Gamasc}
Let $\Gamma$ be a component of a translation quiver which is simply connected, then
$\Gamma$ is a component with length.
\end{thm}

It is not hard to see that  if  $H$ is a hereditary algebra of  Dynkin type, then
the orbit graph of $\Gamma_{H}$ is  of tree type. In this case, the orbit graph is defined as follows; there is a vertex for each indecomposable projective $H$-module and there is a edge between two different vertices  if there exists an irreducible morphism between the indecomposable projective $H$-modules.
Hence, $\Gamma_{H}$ is simply connected and by Theorem \ref{Gamasc}, we get the result.

\begin{prop} \cite{CG} \label{camiglar}
Let $H$ be a finite dimensional algebra of finite representation type. Then, the paths of irreducible morphisms
from $P_{a}$ to $I_{a}$ and  from $P_{b}$ to $I_{b}$ are of the same length, where $P_{a}, I_{a}, P_{b}, I_{b}$, are the indecomposable projective and injective  $H$-modules, corresponding to the vertices $a$ and  $b$ in $Q_{H}$, respectively.
\end{prop}

\begin{proof}
Let $i,j\in Q_{H}$ be vertices  such that there is an arrow from $i$ to $j$ in $Q_{H}$. Then, there are irreducible morphisms between the indecomposable projective modules  $P_{j}$ and $P_{i}$ and between the indecomposable injective modules $I_{j}$ and $I_{i}$.
Moreover, for each vertex $a\in (Q_{H})_{0}$ there is a non-zero path from $P_{a}$ to $I_{a}$, see \cite[Lemma 2.3]{C}. Then, in $\Gamma_{H}$
we have the following diagram:

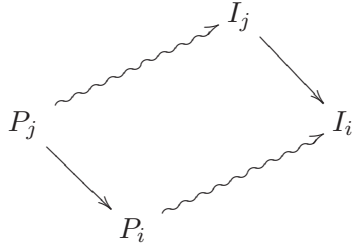
\begin{displaymath}
\xymatrix {& & I_j \ar[rd]  & \\
P_j  \ar[rd]\ar@{~>}[rru] & & & I_i \\
& P_i  \ar@{~>}[rru] & & }
\end{displaymath}

\noindent where by $\rightsquigarrow$ we denote a path of irreducible morphisms.

Now, assume that the path of irreducible morphisms from
$P_{j}$ to $I_{j}$ is of length $r$. Then, the path from
$P_{j} \rightsquigarrow I_{j}\rightarrow I_{i}$ has length $r+1$. Since $\Gamma_H$ is with length then
the path $P_{j}\rightarrow P_{i}\rightsquigarrow I_{i}$
has also length $r+1$. Therefore, the path
$P_{i}\rightsquigarrow I_{i}$ has length $r$.
Since $H$ is connected then all paths of the form
$P_{a}\rightsquigarrow I_{a}$ have the same length.
\end{proof}

Note that if $\mathbf{C_2}({\rm proj} \,H)$ is of finite type then the hereditary algebra
$H$ is of Dynkin type.

\begin{lem}\label{grado de rho y lambda}   Let $H=k Q_H$ with $Q_H$ not semisimple and  $\mathbf{C_2}({\rm proj} \,H)$ is of finite type.
If $i$ and $j$ are a sink  and a source of $Q_H$, respectively then $d_r(\rho_1^{2,i})=d_l(\lambda_1^{2,j})=r_{_H}$, where $r_{_H}$ is the nilpotency index of the radical  of ${\rm mod}\,H$.
\end{lem}

\begin{proof} Since $i$ is a sink of $Q_H$ then, $P_i$ is  a simple projective. Then, there is a non-zero path of irreducible morphisms from $P_i=S_i$ to $I_i$.
Moreover,  from \cite[Teorema 2.6]{C} and \cite{CG} we get that if $r_{_H}$ is the nilpotency index of the radical  of ${\rm mod}\,H$ then the length of a path from $P_i$ to $I_i$ is $r_{_H} -1$.

By \cite[Proposition 3.3]{BSZ:04} there is an equivalence between  ${\rm mod}\,H$ and $\overline{\mathbf{C_2}}({\rm proj} \,H)$ which sends the projective $P_i$ to the complex $(0 \rightarrow P_i)$ in $\mathbf{C_2}({\rm proj} \,H)$ and each module $M$ to its  minimal projective resolution $P_M:P^1_M \rightarrow P^0_M$.

Consider  a  path of irreducible morphisms with non-zero composition $S_i=P_i \rightarrow {M_2} \rightarrow \cdots \rightarrow {M_{m-2}} \rightarrow I_i$ in ${\rm mod} \,H$ of length $r_{_H}-1$. Then, we also have a non-zero path of irreducible morphisms $ (0 \rightarrow P_i) \rightarrow P_{M_2} \rightarrow \cdots \rightarrow P_{M_{m-2}} \rightarrow P_{I_i}$ in $\mathbf{C_2}({\rm proj} \,H)$ of length $r_{_H}-1$.

The almost split sequence in $\mathbf{C_2}({\rm proj} \,H)$ starting in $(0 \rightarrow P_i)$ has more than one indecomposable middle term, because $Q_H$ has more than one vertex. We write it as follows: $(0\rightarrow P_i) \rightarrow (P_i \rightarrow P_i) \oplus P_{N_1} \cdots \oplus P_{N_r} \rightarrow P_M.$
\noindent For some $1 \leq j \leq r$ we have that $P_{N_j}=P_{M_2}$. Moreover, by  \cite[Proposition 6.12]{BSZ:04} we know that $A_2^{-1}(P_{I_i})= (P_i \rightarrow 0)$. Then, there is a path from $(0 \rightarrow P_i)$ to $(P_i \rightarrow 0)$ of length $r_{_H} + 1$. We illustrate the situation as follows:

{\tiny
\begin{displaymath}
\def\objectstyle{\scriptstyle}
\def\labelstyle{\scriptstyle}
\xymatrix {(0 \rightarrow P_i) \ar[r] \ar[rd] &  P_{M_2} \ar[r] &
\cdots \ar[r] &  P_{M_{m-2}} \ar[r] &  P_{I_i} \ar[r] & E \ar[r] &  (P_i \rightarrow 0) \\
 & (P_i \rightarrow P_i)\ar@{~>}[rrrrru]  & & & & & }
\end{displaymath}
}

Any path from  $(0 \rightarrow P_i)$ to $(P_i \rightarrow 0)$  has zero composition, but there is a non-zero morphism from $(P_i \rightarrow P_i)$ to $(P_i \rightarrow 0)$. Since  $\Gamma_{{\mathbf{C_2}({\rm proj} \,H)}}$ is with length then all paths of irreducible morphisms from $(P_i \rightarrow P_i)$ to $(P_i \rightarrow 0)$ have the same length, equal to  $r_{_H}$. Moreover, there is a non-zero path of irreducible morphisms between indecomposable complexes from $(P_i \rightarrow P_i)$ to $(P_i \rightarrow 0)$ in $\Re^{r_{_H}}({\mathbf{C_2}({\rm proj} \,H)}) \backslash \Re^{r_{_H}+1}({\mathbf{C_2}({\rm proj} \,H)})$. The cokernel of $\rho_1^{2,i}:(0 \rightarrow P_i) \rightarrow (P_i \rightarrow P_i)$ is $(P_i \rightarrow 0)$ then  $d_r(\rho_1^{2,i})=r_{_H}$. Dually, we get that $d_l(\lambda_1^{2,j})=r_{_H}$.
\end{proof}

\begin{nota} \label{ell} {\em Let $H=kQ_H$ be a hereditary representation-finite algebra. We denote by $\ell_{j,i}$ the number of arrows
of a path from a source vertex $i$ to a  sink vertex $j$  of $Q_H$ and  $\ell=max\{ \ell_{j,i} \}$.}
\end{nota}

\begin{lem}\label{PiPj} Let $H= k Q_H$ be a representation-finite hereditary algebra. If $f:X\rightarrow Y$ is a non-zero morphism between projective  $H$-modules then there exist morphisms $\alpha:P_i \rightarrow X$ and $\beta:Y \rightarrow P_j$, with $i$ a sink vertex and $j$ a source vertex of
$Q_H$ such that $\beta f \alpha $ is a non-zero monomorphism.
\end{lem}

\begin{proof} Since $f$ is a non-zero morphism then there is an indecomposable projective summand $P$ of $X$ such that, $f(P)\neq 0$.
Consider $\eta:P \rightarrow X$ the inclusion morphism.

On the other hand, there is an indecomposable direct summand  $P'$ of $Y$ such that if $\pi:Y \rightarrow P'$ is the projection morphism then $\pi f(P)\neq 0$.
Therefore, $\pi f \eta \neq 0$. Moreover, $\pi f \eta$ is a monomorphism.

Since  $P$ and $P'$ are indecomposable projective in ${\rm mod}\,H$ there are vertices $t$ and $h$ in  $Q_H$ such that $P=P_t$ and $P'=P_h$.
Moreover, there is a sink $i$,  a source $j$  and a path $\xymatrix {j
\ar@{~>}[r]& h \ar@{~>}[r]& t  \ar@{~>}[r] & i }$ in $Q_H$. Hence, there are non-zero monomorphisms
$\mu:P_i \rightarrow P_t$ and $\nu:P_h \rightarrow P_j$. Let $\alpha= \eta \mu$ and $\beta=\nu \pi$. Then,
 $\beta f \alpha = (\nu \pi)f( \eta \mu) = \nu (\pi f \eta) \mu \neq 0$ and we get that  $\beta f \alpha $ is a monomorphism.
\end{proof}

\begin{lem} \label{paraC2} Let $H= k Q_H$ be a representation-finite hereditary algebra and $f:X \rightarrow Y$  is a non-zero morphism between indecomposable complexes in $\mathbf{C_2}({\rm proj} \,H)$. Then, there exists a sink $i$, a source $j$ in $Q_H$ and morphisms $\psi$ and $\phi$ such that one of the following conditions hold:
\begin{itemize}
\item [a)] $\psi f \phi \in \Re_{\mathbf{C_2}({\rm proj} \,H)}^{r_{_H} + l_{j,i}}  \backslash \Re_{\mathbf{C_2}({\rm proj} \,H)}^{r_{_H} + l_{j,i}+1}$, with $\phi: T(P_i) {\rightarrow} X $ and ${\psi}: Y { \rightarrow} J_1(P_j)$.

  \item [b)] $\psi f \phi \in \Re_{\mathbf{C_2}({\rm proj} \,H)}^{l_{j,i}}  \backslash \Re_{\mathbf{C_2}({\rm proj} \,H)}^{l_{j,i}+1}$, with  $\phi: T(P_i){\rightarrow} X$ and ${\psi}: Y { \rightarrow} T(P_j)$.

  \item [c)] $\psi f \phi \in \Re_{\mathbf{C_2}({\rm proj} \,H)}^{l_{j,i}}  \backslash \Re_{\mathbf{C_2}({\rm proj} \,H)}^{l_{j,i}+1}$, with ${\phi}: S(P_i) {\rightarrow} X $ and ${\psi}: Y {\rightarrow} S(P_j)$.

  \item [d)] $\psi f \phi \in \Re_{\mathbf{C_2}({\rm proj} \,H)}^{r_{_H} + l_{j,i}}  \backslash \Re_{\mathbf{C_2}({\rm proj} \,H)}^{r_{_H} + l_{j,i}+1}$, with $\phi: J_1(P_i) {\rightarrow} X$ and $ {\psi}: Y {\rightarrow} S(P_j)$.
\end{itemize}
\noindent where $r_{_H}$ is the nilpotency index of the radical of ${\rm mod} \,H$  and  $l_{j,i}$ as in  \rm{(\ref{ell})}.
\end{lem}

\begin{proof} Let $f=\{f^1,f^2\}:X \rightarrow Y$  be a non-zero morphism between indecomposable complexes $X=X^1 \stackrel{d_X}{\rightarrow} X^2$ and $Y=Y^1 \stackrel{d_Y}{\rightarrow} Y^2$ in $\mathbf{C_2}({\rm proj} \,H)$. First, consider the case where $f^2 \neq 0$.
By Lemma \ref{PiPj}, there is a sink $i$,  a source $j$, $\alpha:P_i \rightarrow X^2 $ and  $\beta:Y^2 \rightarrow P_j $ such that $\beta f^2 \alpha \neq 0$ and $\beta f^2 \alpha$ is a monomorphism. If $Y^1 \neq 0$ we consider the non-zero composition $\psi f \phi$ in $\mathbf{C_2}({\rm proj} \,H)$, where $\phi=\{0, \alpha\}: T(P_i) \rightarrow X$ and $\psi =\{\beta d_Y, \beta\}:Y \rightarrow J_1(P_j)$. Then, $\lambda_1^{2,j} (\psi f \phi)=0$,  where $\lambda_1^{2,j}: J_1(P_j) \rightarrow S(P_j)$ is the irreducible morphism introduced in $(1.5)$
Hence, $\psi f \phi$ factors through $T(P_j)$ which is the kernel of $\lambda_1^{2,j}$. If $t:T(P_j) \rightarrow J_1(P_j)$ is the kernel morphism of $\lambda_1^{2,j}$  then, there is a morphism $h:T(P_i) \rightarrow T(P_j)$ such that  $\psi f \phi = t h $.  We illustrate the situation with the following diagram:

$$\xymatrix @!0 @R=10mm @C=16mm{T(P_i) \ar[r]^{\phi} \ar[rd]_h & X \ar[r]^f & Y \ar[r]^{\psi} & J_1(P_j) \ar[r]^{\lambda_1^{2,j}} & S(P_j).\\ & T(P_j) \ar[rru]_t & }
$$

\noindent By  Lemma \ref{grado de rho y lambda}, $d_l(\lambda_1^{2,j})=r_{_H}$. Since  $T(P_j)$ is the cokernel of $\lambda_1^{2,j}$, by Remark \ref{conucleo} we have that $t \in \Re_{\mathbf{C_2}({\rm proj} \,H)}^{r_{_H}} (T(P_j), J_1(P_j))\backslash \Re_{\mathbf{C_2}({\rm proj} \,H)}^{r_{_H}+1}(T(P_j), J_1(P_j))$.

On the other hand, there is a non-zero path of irreducible morphisms from $P_i$ to $P_j$ in ${\rm mod}\,H$ of length $\ell_{j,i}$. By \cite[Proposition 3.2]{CPS1} there is a non-zero path of irreducible morphisms from $T(P_i)$ to $T(P_j)$ in $\mathbf{C_2}({\rm proj} \,H)$ of length $\ell_{j,i}.$
Since $H$ is of finite representation type then so is $\mathbf{C_2}({\rm proj} \,H)$. Accordingly, the latter has a generalized standard connected Auslander-Reiten quiver, then $h \in \Re_{\mathbf{C_2}({\rm proj} \,H)}^{\ell_{j,i}} (T(P_i), (T(P_j)) \backslash \Re_{\mathbf{C_2}({\rm proj} \,H)}^{\ell_{j,i}+1}(T(P_i), (T(P_j))$. Furthermore, the composition $\psi f \phi$ is equal to the composition $t h$ where $th \in \Re_{\mathbf{C_2}({\rm proj} \,H)}^{r_{_H}+ \ell_{j,i}} (T(P_i), J_1(P_j)) \backslash \Re_{\mathbf{C_2}({\rm proj} \,H)}^{r_{_H} + \ell_{j,i}+1}(T(P_i), J_1(P_j))$, proving Statement $(a)$.

Consider the case where $f^2 \neq 0$ and $Y^1=0$. Assume that $X^1 \neq 0$. By $(1.4)$ we have that $f^2$ is an irreducible morphism or else is a split epimorphism. Note that these two cases can not occur. In the first case, since $f^2 d_X=0$ we get that $d_X =0$ a contradiction. In the second case, we get that $f$ splits a contradiction.  Then $X^1=0$. Hence, $X=T(P)$ and  $Y=T(P')$. By Lemma \ref{PiPj}, there  are non-zero paths of irreducible morphisms
 $\alpha:P_i \rightarrow P $ and $\beta:P' \rightarrow P_j$ such that $\beta f^2 \alpha \neq 0$. Consider the  non-zero composition $\psi f \phi$ in $\mathbf{C_2}({\rm proj} \,H)$, where $\phi=\{0,\alpha\}:T(P_i) \rightarrow X$ and  $\psi=\{0, \beta\}: Y \rightarrow T(P_j)$.
\noindent Then,  $\psi f \phi \in \Re_{\mathbf{C_2}}^{\ell_{j,i}} (T(P_i), T(P_j)) \backslash \Re_{\mathbf{C_2}}^{\ell_{j,i}+1}(T(P_i), T(P_j))$, where $\ell_{j,i}$ is the length of the non-zero path of irreducible morphisms from $T(P_i)$ to $T(P_j)$ in $\mathbf{C_2}({\rm proj} \,H)$, proving $(b)$.

Consider the case where $f^2 =0$. Then, $f^1 \neq 0$.  If  $X= S(P)$ and $Y=S(Q)$ we choose $P_i$ and $P_j$ as in the above case. By Lemma \ref{PiPj}, there are morphisms   $\alpha:P_i \rightarrow P$ and  $\beta:Q \rightarrow P$ such that $\beta f \alpha$ is a monomorphism.
On the other hand, since there is a non-zero path of irreducible morphisms from  $P_i$ to  $P_j$ in ${\rm mod}\,H$  of length  $\ell_{j,i}$  then
 $\{\beta,0\} f \{\alpha,0\} \in \Re_{\mathbf{C_2}({\rm proj} \,H)}^{\ell_{j,i}} (S(P_i), S(P_j))  \backslash \Re_{\mathbf{C_2}({\rm proj} \,H)}^{\ell_{j,i}+1}(S(P_i), S(P_j))$, where  $\{\alpha,0\}:  S(P_i)\rightarrow S(P)$ and  $\{\beta,0\}: S(Q) \rightarrow S(P_j)$ in $\mathbf{C_2}({\rm proj} \,H)$,  proving $(c)$.

If   $X^2\neq 0$,  $Y= S(P')$, again by  Lemma \ref{PiPj} there is a sink $i$, a source $j$,
$\alpha:P_i \rightarrow X^1 $ and $\beta:P' \rightarrow P_j $ such that $0 \neq \beta f^1 \alpha$ is a monomorphism. If $\phi=\{ \alpha, d_X \alpha\}:J_1(P_i) \rightarrow X$ and $\psi= \{\beta,0\}:Y \rightarrow S(P_j) $ then $\psi f \phi \neq 0$.

\noindent Consider the irreducible morphism  $\rho_1^{2,i}: T(P_i) \rightarrow J_1(P_i)$. Then, $(\psi f \phi) \rho_1^{2,i}=0$. Hence, $\psi f \phi $ factors through $S(P_i)$ which is the  cokernel of $\rho_1^{2,i}$. If $t':J_1(P_i) \rightarrow S(P_i)$ is the cokernel morphism then there exists $h':S(P_i) \rightarrow S(P_j)$ such that $h' t' = \psi f \phi$. We illustrate the situation with the following diagram:

{\tiny
$$\xymatrix @!0 @R=10mm @C=16mm{T(P_i) \ar[r]^{\rho_1^{2,i}} & J_1(P_i) \ar[r]^{\phi} \ar[rrd]_{t'} & X \ar[r]^f & Y \ar[r]^{\psi} & S(P_j). \\ & & & S(P_i) \ar[ru]_{h'} &
}
$$ }

\noindent Then, $h' \in \Re_{\mathbf{C_2}({\rm proj} \,H)}^{\ell_{j,i}} (S(P_i), S(P_j)) \backslash \Re_{\mathbf{C_2}({\rm proj} \,H)}^{\ell_{j,i}+1}(S(P_i), S(P_j))$. Since $d_r(\rho_1^{2,i})=r_{_H}$, then by Remark \ref{conucleo} we have that $t' \in \Re_{\mathbf{C_2}({\rm proj} \,H)}^{r_{_H}} (J_1(P_i), S(P_i)) \backslash \Re_{\mathbf{C_2}({\rm proj} \,H)}^{r_{_H} +1}(J_1(P_i), S(P_i))$. Therefore, the composition  $\psi f \phi$ is equal to the composition $h' t'$ where $$h' t' \in \Re_{\mathbf{C_2}({\rm proj} \,H)}^{r_{_H} + \ell_{j,i}} (J_1(P_i), S(P_j)) \backslash \Re_{\mathbf{C_2}({\rm proj} \,H)}^{r_{_H} + \ell_{j,i}+1}(J_1(P_i), S(P_j))$$
\noindent  proving $(d)$.
\end{proof}

Similar results to Lemma \ref{grado de rho y lambda} and Lemma \ref{paraC2} hold in a category $\mathbf{C_3}({\rm proj} \,H)$ of finite type.
These properties allow us to compute the nilpotency index of the radical of this  category.

\begin{lem}\label{gradoderho y lambdaC3}  Let $H=k Q_H$, where $Q_H$ is not semisimple. Let $\mathbf{C_3}({\rm proj} \,{ H})$ be of finite type. If $i$ is a sink and $j$ is a source in  $Q_H$ then $d_r(\rho_1^{3,i})=d_r(\rho_2^{3,i})=d_l(\lambda_1^{3,j})=d_l(\lambda_2^{3,j})=r_H$, where $r_H$ is the nilpotency index of the radical of ${\rm mod}\,{ H}$.
\end{lem}

\begin{proof} From (\ref{1.1}) we know the shape of the  Auslander-Reiten quiver of $\mathbf{C_3}({\rm proj} \,{ H})$. That is, the quiver is obtain as one copy of the quiver of $\mathbf{C_2}({\rm proj} \,{ H})[0]$ and one copie of the quiver of $\mathbf{C_2}({\rm proj} \,{ H})[-1]$, where $$\mbox{ind}\, {\mathbf{C_2}({\rm proj} \,H )[0]} \cap \mbox{ind}\, {\mathbf{C_2}({\rm proj} \,H )[-1]} = \{ (0 \rightarrow P \rightarrow 0) \mid P \in \mbox{ind}({\rm proj} \,H)\}.$$

 The  embedding $\varphi_0$ and $\varphi_1$ defined in (\ref{1.1})  satisfy that $\varphi_0(\rho_1^{2,i})=\rho_1^{3,i}$ and  $\varphi_1(\rho_1^{2,i})=\rho_2^{3,i}$.  Applying the embedding $\varphi_0$ and  $\varphi_1$, respectively in the proof of Lemma \ref{grado de rho y lambda} we obtain that $d_r(\rho_1^{3,i})=d_r(\rho_2^{3,i})=r_H$. Similarly for the left degree of  $\lambda_1^{3,j}$ and $\lambda_2^{3,j}$.
\end{proof}

\begin{lem}\label{paraC3} Let $\mathbf{C_3}({\rm proj} \,{ H})$ be a category of finite type  and  $f:X \rightarrow Y$ be a non-zero morphism between indecomposable complexes in
$\mathbf{C_3}({\rm proj} \,{H})$. Then, there is a sink $i$,  a source $j$ in $Q_H$ and morphisms $\psi$ and $\phi$ such that one of the following conditions hold.
\begin{itemize}
  \item [a)] $\psi f \phi \in \Re_{\mathbf{C_3}}^{l_{j,i}}  \backslash \Re_{\mathbf{C_3}}^{l_{j,i}+1}$, with  $\phi: S_k(P_i) {\rightarrow} X $, and ${\psi}: Y { \rightarrow} S_k(P_j)$, for some $k\in \{1,2,3\}$.

  \item [b)] $\psi f \phi \in \Re_{\mathbf{C_3}}^{l_{j,i} + r_H}  \backslash \Re_{\mathbf{C_3}}^{l_{j,i} + r_H + 1}$, with  $\phi: T(P_i){\rightarrow} X$ and ${\psi}: Y { \rightarrow} J_2(P_j)$.

  \item [c)] $\psi f \phi \in \Re_{\mathbf{C_3}}^{l_{j,i}+2r_H}  \backslash \Re_{\mathbf{C_3}}^{l_{j,i}+2r_H+1}$, with  ${\phi}: J_2(P_i) {\rightarrow} X $ and  ${\psi}: Y {\rightarrow} J_1(P_j)$.

  \item [d)] $\psi f \phi \in \Re_{\mathbf{C_3}}^{r_H + l_{j,i}}  \backslash \Re_{\mathbf{C_3}}^{r_H + l_{j,i}+1}$, with  $\phi: J_1(P_i) {\rightarrow} X$ and  ${\psi}: Y {\rightarrow} S(P_j)$.
\end{itemize}

\noindent where $r_H$ is the nilpotency index of the radical of ${\rm mod} \,{ H}$  and $l_{j,i}$ defined as in (\ref{ell}).
\end{lem}

 \begin{proof} Let $f=\{f^1,f^2,f^3\}: X \rightarrow Y$ be a morphism between indecomposable complexes in $\mathbf{C_3}({\rm proj} \,{ H})$. Consider three cases,  $(1)$.  $f^1=0$ and $f^2 \neq 0$, $(2)$. $f^1=f^2=0$ and  $f^3 \neq 0$ and $(3)$. $f^1 \neq 0$.

Let $f^1=0$ and $f^2 \neq 0$. By Lemma \ref{PiPj} there are indecomposable projective modules $P_i$ and $P_j$ in ${\rm mod}\,{ H}$, with $i$ a sink, $j$ a source in $Q_H$ and morphisms $\alpha: P_i \rightarrow X^2$ and $ \beta: Y^2 \rightarrow P_j$ such that $\beta f^2 \alpha$ is a non-zero monomorphism.

If $X^1=Y^1=X^3=Y^3=0$, we define $\phi=\{0,\alpha,0\}:S_2(P_i) \rightarrow X$ and  $\psi=\{0,\beta,0\}:Y \rightarrow S_2(P_j)$. Then, $\psi f \phi \in   \Re_{\mathbf{C_3}}^{l_{j,i}}  \backslash \Re_{\mathbf{C_3}}^{l_{j,i}+1}$. Hence, we get that $a)$ holds, for $k=2$.

Otherwise, we consider the morphisms $\phi=\{0,\alpha,d_{X^2} \alpha\}: J_2(P_i) \rightarrow X$ and $\psi=\{\beta d^1_Y, \beta, 0\}:Y \rightarrow J_1(P_j)$. Then, $\psi f \phi: J_2(P_i) \rightarrow J_1(P_j)$ is non-zero. By Lemma \ref{gradoderho y lambdaC3} we have that $d_r(\rho^{3,i}_2)=r_H$. Therefore, the cokernel morphism $t:J_2(P_i) \rightarrow S_2(P_i)$ satisfies that$t \in \Re_{\mathbf{C_3}}^{r_H}  \backslash \Re_{\mathbf{C_3}}^{r_{H}+1}$. Moreover,
$(\psi f \phi) \rho^{3,i}_2=0$. Then, the morphism $\psi f \phi$  factors through ${\rm Coker}\,\rho^{3,i}_2=(P_i \rightarrow 0\rightarrow 0) =S(P_i)$, that is, there is a morphism $h:S(P_i) \rightarrow J_1(P_j)$ such that $h t = \psi f \phi$.

$$
\def\objectstyle{\scriptstyle}
\def\labelstyle{\scriptstyle}
\xymatrix @!0 @R=10mm @C=16mm{T(P_i) \ar[r]^{\rho^{3,i}_2} & J_2(P_i) \ar[r]^{\phi} \ar[rrd]_{t} & X \ar[r]^f & Y \ar[r]^{\psi} & J_1(P_j)\ar[r]^{\lambda^{3,j}_1} & S_1(P_j) \\ & & & S_2(P_i) \ar@{-->}[ru]_{h} \ar@{-->}[rr]_{h'} & & S_2(P_j) \ar[ul]^{t'}
}
$$

On the other hand, consider the irreducible morphism  $\lambda_1^{3,j}:J_1(P_j) \rightarrow S_1(P_j)$. By Lemma \ref{gradoderho y lambdaC3} we have that  $d_l(\lambda^{3,j}_1)= r_H$. Hence, the kernel morphism  $t':S_2(P_j) \rightarrow J_1(P_j)$ satisfies that $t' \in \Re_{\mathbf{C_3}}^{r_H}  \backslash \Re_{\mathbf{C_3}}^{r_H+1}$. Since $\lambda_1^{3,j} h =0$ then the morphism $h$ factors through the kernel of $\lambda_1^{3,j}$. That is, there is a morphism  $h': S_2(P_i) \rightarrow S_2(P_j) \in \Re_{\mathbf{C_3}}^{l_{j,i}}  \backslash \Re_{\mathbf{C_3}}^{l_{j,i}+1}$ such that $h=t'h'$. Then, $\psi f \phi = h t$ and $h t = t' h' t$ where $ t' h' t \in \Re_{\mathbf{C_3}}^{l_{j,i}+2 r_H}  \backslash \Re_{\mathbf{C_3}}^{l_{j,i}+2 r_H+1}$, since $\Gamma_{\mathbf{C_3}}$ is a generalized standard component with length. Therefore, this proves that $f$ satisfies $c)$.

The proof of the remaining cases  follow similarly.
\end{proof}

\begin{coro}\label{f en radical en C3} Let $H= k Q_H$ be a representation-finite hereditary algebra, $r_{_H}$ the nilpotency index of the radical of
${\rm mod}\,H$ and $\ell$ defined as in \ref{ell}. Let $f:X \rightarrow Y$ is a non-zero morphism  with $X,Y$ indecomposable in $\mathbf{C_n}({\rm proj} \,H)$. The following conditions hold.

\begin{enumerate}
\item[(a)] If  $n=2$ then $f \in \Re_{\mathbf{C_2}({\rm proj} \,H)}^s (X, Y) \backslash \Re_{\mathbf{C_2}({\rm proj} \,H)}^{s+1}(X, Y)$ with $s \leq r_{_H} + \ell$.

\item[(b)] If  $n=3$ then $f \in \Re_{\mathbf{C_3}({\rm proj} \,H)}^s (X, Y) \backslash \Re_{\mathbf{C_3}({\rm proj} \,H)}^{s+1}(X, Y)$ with $s \leq 2r_{_H} + \ell$.
\end{enumerate}
\end{coro}

In \cite{CG}, the author determined the nilpotency index for the radical of ${\rm mod}\,H$, where $H$ is of Dynkin type. The bound was given in terms of the number of vertices of its ordinary quiver. Moreover, in \cite[Theorem 4.10]{Z} D. Zacharia also determined such bounds.

More precisely, we have the following result:

\begin{thm} \label{ChaioGuazzelli}  Let $H= k Q$  be a representation-finite hereditary algebra. If $r_{_H}$ is the nilpotency index of the radical of then:
\begin{itemize}
\item [(a)] If  $Q= A_n$ then $r_{_H} = n$.
  \item [(b)] If $Q= D_n$ then  $r_{_H}= 2n - 3$, for $n \geq 4$.
 \item [(c)] If  $Q= E_n$ for $n=6,7,8$ then $r_{_H}= 11, 17, 29$, respectively.
\end{itemize}
\end{thm}

Now, we are in a position to prove the main result of this section,  Theorem D.

\begin{thm}\label{indice} Let $H$ be a hereditary algebra, and $r_{_H}$ be the nilpotency index of the radical of ${\rm mod}\,H$ and $\ell$ defined as in \ref{ell}. Let $\mathbf{C_n}({\rm proj} \,H)$ be a representation-finite category. Then, the nilpotency index $r$ of the  radical of  $\mathbf{C_n}({\rm proj} \,H)$ is the following:
\begin{enumerate}
\item [(a)] If   $n=2$    then $r=r_{_H} + \ell + 1$.
\item [(b)] If $n \geq 3$ then  $r= 2r_{_H} + \ell + 1$.
\end{enumerate}
\end{thm}

\begin{proof} $(a)$ By Corollary \ref{f en radical en C3} $(a)$, if $f:X \rightarrow Y$ is a non-zero morphism between indecomposable complexes  in $\mathbf{C_2}({\rm proj} \,H)$  then $f \in \Re_{\mathbf{C_2}({\rm proj} \,H)}^s (X,Y) \backslash \Re_{\mathbf{C_2}({\rm proj} \,H)}^{s+1}(X,Y)$, with $s \leq r_H + \ell$, where $r_H$ is the nilpotency index of the radical of  ${\rm mod}\,H$ and $\ell$ defined as in \ref{ell}.
We claim that $r= r_{_H} + \ell + 1$.
We  know that $\Re ^{r_{_H} + \ell +1} ({\mathbf{C_2}({\rm proj} \,H)})=0$. Let prove that $\Re ^{r_{_H} + \ell} ({\mathbf{C_2}({\rm proj} \,H)})\neq 0$. Consider a sink $i$ and a source $j$ in $Q_H$ such that $\ell = \ell_{j,i}$. Since there is a path from  $j$ to $i$ in $Q_H$ of length $\ell$ then there is a morphism $d: S(P_i)\rightarrow S(P_j)$ in $\mathbf{C_2}({\rm proj} \,H)$ by \cite[Proposition 3.4]{CPS1}. Moreover, since by Lemma \ref{grado de rho y lambda} $d_r(\rho_1^{2,i})=r_{_H}$, and by  Remark \ref{conucleo} the cokernel morphism $v:J_1(P_i) \rightarrow S(P_i)$ satisfies that $v \in \Re_{\mathbf{C_2}({\rm proj} \,H)}^{r_{_H}} (J_1(P_i), S(P_i)) \backslash \Re_{\mathbf{C_2}({\rm proj} \,H)}^{r_{_H} +1}(J_1(P_i), S(P_i))$. Then, $$d v \in \Re_{\mathbf{C_2}({\rm proj} \,H)}^{r_{_H} + \ell} (J_1(P_i), S(P_j))\backslash \Re_{\mathbf{C_2}({\rm proj} \,H)}^{r_{_H} + \ell +1}(J_1(P_i), S(P_j))$$ proving the result.

$(b)$
By Corollary \ref{f en radical en C3} $(b)$, if $f:X \rightarrow Y$ is a morphism between indecomposable in $\mathbf{C_3}({\rm proj} \,H)$ then $f \in \Re_{\mathbf{C_3}({\rm proj} \,H)}^s (X, Y) \backslash \Re_{\mathbf{C_3}({\rm proj} \,H)}^{s+1}(X, Y)$, with $s \leq 2 r_{_H} + \ell$.  We claim that $r =2 r_{_H} + \ell +1$.  Consider a sink $i$ and a source $j$ in $Q_H$ such that $\ell = \ell_{j,i}$. Then, there is a morphism $k \in \Re_{\mathbf{C_3}({\rm proj} \,H)}^{\ell} (S_2(P_i), S_2(P_j))\backslash \Re_{\mathbf{C_3}({\rm proj} \,H)}^ {\ell +1}(S_2(P_i), S_2(P_j))$.
Moreover, $ d_r(\rho_2^{3,i})=r_{_H}$, where $\rho_2^{3,i}: R_2(P_i) \rightarrow J_2(P_i)$. Then,   $t \in \Re_{\mathbf{C_3}({\rm proj} \,H)}^{r_{_H}} (J_2(P_i), S_2(P_i)) \backslash \Re_{\mathbf{C_3}({\rm proj} \,H)}^{r_{_H}+1}(J_2(P_i), S_2(P_i))$, where $t$ is the  cokernel morphism of $\rho_2^{3,i}$.
Moreover,   $d_l(\lambda_1^{3,j})=r_{_H}$ with  $\lambda_1^{3,j}: J_1(P_j) \rightarrow L_1(P_j)$. Therefore,  $t' \in \Re_{\mathbf{C_3}({\rm proj} \,H)}^{r_{_H}} (S_2(P_j), J_1(P_j)) \backslash \Re_{\mathbf{C_3}({\rm proj} \,H)}^ {r_{_H} +1}(S_2(P_j), J_1(P_j)) $ where $t'$ is the kernel morphism of $\lambda^{3, j}_1$. Since $t' k t \neq 0$ and  $\Gamma_{C_3({\rm proj} \,H)}$ is generalized standard with length then the composition $t' k t : J_2(P_i) \rightarrow  J_1(P_j) $ is such that
$t' k t \in \Re_{\mathbf{C_3}({\rm proj} \,H)}^{2 r_{_H} + \ell} (J_2(P_i),  J_1(P_j)) \backslash \Re_{\mathbf{C_3}({\rm proj} \,H)}^{ 2r_{_H} + \ell +1} (J_2(P_i),  J_1(P_j))$. We illustrate the situation as follows:

{\tiny
$$
\xymatrix @!0 @R=10mm @C=16mm {J_2(P_i) \ar[d]^{t} & 0 \ar[r] \ar[d] & P_i \ar[r]^{1} \ar[d]^{1} & P_i \ar[d]^{}
\\S_2(P_i) \ar[d]^{k} & 0 \ar[r]\ar[d] & P_i \ar[d]\ar[r] &  0\ar[d]
\\S_2(P_j) \ar[d]^{t'} & 0 \ar[r]\ar[d] & P_j \ar[d]^{1} \ar[r] &  0 \ar[d]
\\J_1(P_j) & P_j \ar[r]^{1} &  P_j \ar[r] & 0}
$$
}

\noindent Hence, the above path is the longest, therefore $\Re^{ 2r_{_H} + \ell} ({\mathbf{C_3}({\rm proj} \,H)})\neq 0$, and the  nilpotency index of $\Re ({\mathbf{C_3}({\rm proj} \,H)})$  is $2r_{_H} + \ell +1$.

Finally, for each $j$ and $k \in \{0, \cdots, n-2 \}$, there are no non-zero morphisms from objects in $\mathbf{C_2}({\rm proj} \,H)[-j]$ to objects in $\mathbf{C_2}({\rm proj} \,H)[-k]$ unless $k=j-1$, and by $(\ref{1.1} )$ we have that

$${\rm ind}(\mathbf{C_n}({\rm proj} \,H))=  \cup^{n-2}_{j=0}{\rm
ind}(\mathbf{C_2}({\rm proj} \,H)[-j])$$
\noindent  then the  nilpotency index of the radical of $ \mathbf{C_n}({\rm proj} \,H)$ is the same that the  nilpotency index of the radical of $\mathbf{C_3}({\rm proj} \,H)$.\end{proof}

As an immediate consequence we state the following corollary.

\begin{coro} Let $\mathbf{C_n}({\rm proj} \,H)$ be a category of finite type with $H$  a hereditary algebra. Consider $r=r_{_H} + \ell + 1$ if $n=2$ and $r= 2r_{_H} + \ell + 1$ if $n \geq 3$.
Then, the composition of $r+1$ irreducible morphisms between complexes in $\mathbf{C_n}({\rm proj} \,H)$ is zero.
\end{coro}

\end{document}